\newcommand{\iMaslov}{\mu_{\scriptscriptstyle{\mathrm{Maslov}}}}

\newcommand{\irel}{\mu_{\scriptscriptstyle{\mathrm{rel}}}}
\newcommand{\ispec}{\mu_{\scriptscriptstyle{\mathrm{spec}}}}

\newcommand{\iRcon}{i_{\scriptscriptstyle{\mathrm{con}}}}
\newcommand{\iMor}{\mu_{\scriptscriptstyle{\mathrm{Morse}}}}
\newcommand{\icon}{\mu_{\scriptscriptstyle{\mathrm{con}}}}
\newcommand{\isf}{\mu_{\scriptscriptstyle{\mathrm{spec}}}}

\newcommand{\spfl}{\mathrm{sf}}
\newcommand{\Imm}{\mathrm{Im}}

\newcommand{\Ker}{\mathrm{ker}}

\newcommand{\sgn}{\mathrm{sign}}

\newcommand{\Id}{\mathrm{Id}}

\newcommand{\Dds}{\tfrac{D}{ dx}}

\newcommand{\Ddt}{\tfrac{ D}{dx}}
\newcommand{\Ddtt}{\tfrac{ D^2}{ dx^2}}

\newcommand{\R}{\mathbb R}
\newcommand{\Z}{\mathbb Z}
\newcommand{\C}{\mathbb C}

\newcommand{\RR}{{\mathbb R}}

\newcommand{\be}{\begin{equation}}
\newcommand{\ee}{\end{equation}}
\newcommand{\bea}{\begin{eqnarray}}
\newcommand{\eea}{\end{eqnarray}}

\newcommand{\equ}[1]{(\ref{#1})}

\newcommand{\ve}{\varepsilon}
\def\virgA{``}

\documentclass[11pt,final, a4paper]{amsart}
\usepackage{amssymb}

\title[A Morse Index Theorem]{A Morse Index
Theorem for Perturbed \\Geodesics on semi-Riemannian
Manifolds}

\author[M. Musso]{Monica Musso }
\address{Dipartimento di Matematica
\hfill\break\indent Politecnico di Torino \hfill\break\indent
Torino, To, Italy} \email{musso@calvino.polito.it}
\author[J. Pejsachowicz]{Jacobo Pejsachowicz}
\address{Dipartimento di Matematica
\hfill\break\indent Politecnico di Torino \hfill\break\indent
Torino, To, Italy} \email{jacobo@polito.it}
\author[A. Portaluri]{Alessandro Portaluri}
\address{Departamento de Matem\'atica,\hfill\break\indent
Instituto de Matem\'atica e Estat\'\i stica\hfill\break\indent
Universidade de S\~ao Paulo, \hfill\break\indent Rua do Mat\~ao
1010, CEP 05508-900, S\~ao Paulo, SP\hfill\break\indent Brazil}
\email{portalur@ime.usp.br} \curraddr{Dipartimento di Matematica,
Politecnico di Torino, Italy.}


\numberwithin{equation}{section}

\begin{document}

\theoremstyle{plain}\newtheorem{teo}{Theorem}[section]
\theoremstyle{plain}\newtheorem{proposition}[teo]{Proposition}
\theoremstyle{plain}\newtheorem{lemma}[teo]{Lemma}
\theoremstyle{plain}\newtheorem{cor}[teo]{Corollary}
\theoremstyle{definition}\newtheorem{definition}[teo]{Definition}
\theoremstyle{remark}\newtheorem{rem}[teo]{Remark}
\theoremstyle{definition}\newtheorem{example}[teo]{Example}
\theoremstyle{plain}\newtheorem*{convention}{Convention}
\theoremstyle{definition}\newtheorem*{defin0}{Definition}

\maketitle

\begin{abstract}

Perturbed geodesics are trajectories of particles moving  on a
semi-Riemannian  manifold in the presence of a potential. Our
purpose here is to extend to perturbed geodesics on
semi-Riemannian manifolds the well known Morse Index Theorem. When
the metric is indefinite, the  Morse index of the energy
functional becomes  infinite and hence, in order to obtain a
meaningful statement, we substitute the Morse index by its
relative form, given by the spectral flow of an associated family
of index forms.  We also introduce a new counting  for conjugate
points, which need not to be isolated in this context, and  prove
that  our generalized  Morse index equals the total number of
conjugate points.  Finally we  study the relation with the Maslov
index of the flow induced on the Lagrangian Grassmannian.
\end{abstract}

\bigskip
\bigskip
     \centerline{Dedicated to Kazimierz G\c eba}
\bigskip
\bigskip

\section{Introduction}

A semi-Riemannian manifold  is a smooth $n$-dimensional manifold
$M$ endowed with a (pseudo) metric given by a nondegenerate
symmetric two-form $g$ of constant index $\nu$.  We denote by  $D$
the associated Levi-Civita connection and by $\Ddt$ the covariant
derivative of a vector field along a smooth curve $\gamma$. Let
$I$ be an interval on the real line. Let    $V$ be  a smooth
function defined  on $I\times M$. A {\em perturbed geodesic\/}
abbreviated as {\em p-geodesic\/} is a smooth curve $\gamma\colon
I \rightarrow M$ which  satisfies the differential equation
\begin{equation}\label{2.2}
\Ddt \gamma '(x) +\nabla V(x,\gamma(x))=0
\end{equation}
 where  $\nabla V$ denotes
gradient of $V (x,-) $ with respect to the metric $g.$

 From the viewpoint  of analytical dynamics, the data $(g,V)$
define a mechanical system on the manifold $M$, with kinetic
energy $\frac{1}{2}g(v ,v)$ and time dependent  potential energy
$V.$ Solutions of the differential equation \equ{2.2} are
trajectories of particles moving on the semi-Riemannian manifold
in the presence of the potential $V$. If the potential vanishes we
get trajectories of  free  particles and hence geodesics on $M$.
This motivates the  suggestive name, \virgA perturbed geodesics",
already adopted in \cite{W} for the periodic case. If the
potential $V$ is time independent then, modulo reparametrization,
perturbed geodesics become geodesics of the Jacobi metric
associated to $(g,V).$
The total energy $$ e = \frac{1}{2}g
(\gamma (x)) ( \gamma '(x) , \gamma '(x) ) + V (\gamma(x) )$$ is
constant along such  trajectories and when $V$ is bounded from
above the solutions of \equ{2.2} with energy $e$ greater than
$\sup_{m \in M} V(m)$ are nothing but reparametrized geodesics for
metric $ [ e-V]g$ on $M$ with total energy one \cite{AM}.

In what follows we will consider perturbed geodesics connecting
two given points of $M$ and will normalize the domain  by taking
as $I$ the interval $[0,1]$.

A vector field $\xi$ along $\gamma$  is called a {\em Jacobi
field\/} if it verifies the  linear differential equation

\begin{equation}\label{2.5}
\Ddtt\xi (x) + R(\gamma ' (x) , \xi (x)) \gamma '(x) + D_{\xi
(x)}\nabla V(x, \gamma(x))=0,
\end{equation}
where $R$ is the curvature tensor of $D.$

Given a $p$-geodesic  $\gamma$, an instant $x \in (0,1]$ is said
to be a {\em conjugate instant\/} if there exists at least one non
zero Jacobi field with  $\xi(0)=\xi (x)=0.$
    The corresponding  point $ q=\gamma(x)$ on $M$ is said to be a
{\em conjugate point\/} to the point $p =\gamma(0)$ along
$\gamma.$

Let  ${\mathcal I}$ be the $n$-dimensional vector space of all
Jacobi fields along $\gamma$ verifying $\xi(0)=0$.  The number $
m(x ) = \dim\, \{ \xi \in {\mathcal I}: \xi (x ) =0 \}$ is called
the {\em geometric multiplicity\/} of $x$.  Thus  $x$ is a
conjugate instant if and only if $m(x) >0$. Let $\mathcal I [x] =
\{ \xi (x): \xi \in \mathcal I \} \subset T_{\gamma (x)} M$.
Denoting with $\perp$ the orthogonal with respect to the metric
$g$, the rank theorem applied to the homomorphism $ \xi \in
{\mathcal I}\, \rightarrow \xi(x)\, $ gives $m(x) = {\rm codim} \,
\mathcal I [x] = \dim \mathcal I[x]^{\perp}.$

We will say that a  conjugate instant $x$ is {\em regular}
  if the restriction of the form $g_{\gamma(x) } $ to ${\mathcal
I}[x]^{\perp}$ is  a non degenerate  quadratic form and that
$\gamma $ is a {\em regular\/} $p$-geodesic if all the conjugate
instants along $\gamma $ are regular. It is easy to see that
regular conjugate instants are isolated and therefore any regular
$p$-geodesic $\gamma$ has only a finite number of conjugate
instants.

The {\em conjugate index\/} of a regular $p$-geodesic  $\gamma$ is
defined by
\begin{equation} \label{icon}
\iRcon(\gamma ) =  \sum_{x \in (0,1]} \sgn\,\big( g|_{{\mathcal
I}[x]^{\bot}}\big)
\end{equation}
where $\sgn$  denotes the signature of a quadratic form.

It is easy to see  that  both the regularity and $ \iRcon(\gamma
)$ are preserved under small enough $p$-geodesic perturbations of
a regular $p$-geodesic. For geodesics  on  a semi-Riemannian
manifold this was proved in \cite{H} and \cite{GMPT}, where it was
also shown that the sum  in \equ{icon} is not a topological
invariant of $\gamma$ if the geodesic is not regular.

In what follows we will  shortly describe some  known results in
the geodesic case. If the metric is Riemannian the conjugate
points along a  geodesic are always regular and   $\sgn\, \big(
g|_{{\mathcal I}[x]^{\bot}}\big)=m(x).$
  The regularity of all conjugate points implies then that
$\iRcon(\gamma )$ is a topological invariant. The celebrated Morse
index theorem states that  the Morse index $\iMor(\gamma)$ of the
geodesic $\gamma$, considered  as a critical point of the energy
functional, equals the total  multiplicity  $ \sum_{x \in ]0,1]}
m(x)$ of conjugate points along the geodesic.

The regularity  property and  the Morse index theorem,
appropriately restated, continue to be true for time-like and
light-like geodesics on Lorentz manifolds. Everything turns to be
as before, eventually after passing to a quotient  space of the
space of Jacobi fields (see \cite{B} for the proof of this result
and its applications  to relativity). However the above approach
breaks down for space-like geodesics or for geodesics of any
causal character on  general semi-Riemannian  manifolds. In this
case not only the Morse index of a geodesic fails to be finite but
also the conjugate points  can accumulate. Moreover they can
disappear under a small perturbation. Thus, in order to formulate
the Morse index theorem for geodesics and  $p$-geodesics on
semi-Riemannian manifolds, one has to extend \equ{icon} to a
topological invariant of general $p$-geodesics and find some kind
of renormalized Morse index as a substitute for the right hand side.

The first breakthrough for this problem  in the geodesic case was
obtained by  Helfer in the paper \cite{H}. His substitute  for the
Morse index  is a spectral index defined as the sum of Krein
signatures of negative real eigenvalues of the Jacobi differential
operator, viewed as a self-adjoint operator on a Krein space.
Recently, in a series of papers, Piccione, Tausk and their
collaborators \cite{MPT,PT, GMPT, CPPT, PPT}  initiated a sort of
critical review of Helfer's work, improving and completing his
results. They fixed a technical gap in Helfer's argument in
\cite{CPPT}.
  In \cite{MPT}, following Helfer, they defined conjugate index by
equation \equ{icon} for any geodesic, irrespective of whether it
is regular or not, but  they gave an example which shows that this
naive definition  does not work in the non regular case. However
the correct expression for the conjugate index of a degenerate
conjugate point was not given,
  and therefore  they conclude that the Morse index theorem holds true
for regular geodesics only.  Their main  result  in \cite{PT,
GMPT} is an expression for the spectral index as the difference
between Morse indices of the index form (the Hessian) restricted
to some special subspaces of the domain. Similar   decomposition
was  already stated by Helfer in \cite{H} without proof.

What we  propose here  is a  different definition  of the two
sides in the Morse  index theorem which extends to perturbed
geodesics as well.
  This allows us to handle  both the problem with the Morse index and
the accumulation of conjugate points, providing a new and, we
believe, interesting proof of this theorem.

The Hessian $h_\gamma $ of the energy functional  at a
$p$-geodesic $\gamma$  is a bounded Fredholm  quadratic form. If
moreover $h_\gamma$ is nondegenerate, the $p$-geodesic  $\gamma$
is called nondegenerate. Together with a nondegenerate
$p$-geodesic $\gamma$ we consider a path $\tilde\gamma$   of
perturbed geodesics canonically  induced by $\gamma$ on the
manifold $\Omega(M)$ of all $H^1$-paths on $M$ (see Section 3). As
generalized Morse index of $\gamma$ we take the negative of the
spectral flow of  the family of Hessians of  the energy functional
along the canonical path $\tilde\gamma$.

Roughly speaking, the spectral flow of a path of  Fredholm
quadratic forms or, what is the same, the spectral flow of the path
 $\{A_t \} _{t\in [a,b]} $ of self-adjoint Fredholm  operators  arising in the Riesz
representation of the forms, is the integer given by the number
of negative eigenvalues of $A_a$ that become positive as the
parameter $t$  goes from $a$ to $b$ minus the number of positive
eigenvalues of $A_a$ that   become negative. It is easy to see
that if one  of (and hence all) the operators in the path have a
finite Morse index, then the spectral flow of  a path $A$ is
nothing but the difference between the Morse indices at the end
points. Thus spectral flow appears to be the right  substitute of
the Morse index  in the framework of strongly  indefinite
functionals. For example, in \cite{F},  it substitutes the Morse
index  in defining a grading for the Floer homology groups. In
\cite{P} it plays the same role  in the formulation of a
bifurcation theorem for critical points for strongly indefinite
functionals. In general the spectral flow depends on the homotopy
class of the whole path and  not only on its end-points. However
in the specific case of the energy functional associated to  a
mechanical system things are simpler. In this framework, spectral
flow depends only on the endpoints of the path and therefore  it
can be considered as a relative form of Morse index in the sense
of \cite{FF}.

We also change the  counting of conjugate points. The papers
\cite{H,PT} follow the lines of the proof of the Morse index
theorem by Duistermaat \cite{D},  which essentially  leads  to
consider as the \virgA total number of  conjugate points" along
$\gamma$ the Maslov index of the flow line  induced  on the
Lagrangian Grassmanian by the associated Hamiltonian flow. We
define  the conjugate  index $\icon$  by means of a suspension of
the complexified family of boundary value problems defining Jacobi
fields. The resulting  boundary value problem is parameterized by
points of the complex plane. The conjugate instants are in one to
one correspondence with points in the complex plane where the
determinant of the associated fundamental matrix vanishes. This
determinant defines a smooth map from an open subset of the
complex plane into the plane. A well known topological invariant
that counts algebraically the zeroes of a map is the Brouwer
degree. In order to keep signs according to  the Riemannian case
we define the conjugate index $\icon$ as the minus the degree of
this map. With  this said, our main theorem  takes the standard
form.

\begin{teo}\label{TEOREMA}
Let $(M , g)$ be a semi-Riemannian manifold, $V\colon I \times
M\to \R$ a smooth potential  and $\gamma \colon I\rightarrow M
$ a nondegenerate perturbed geodesic.  Then $ \isf(\gamma) \, =
\icon (\gamma ).$
\end{teo}

While  $\ispec(\gamma)$ has an intrinsic, i.e. coordinate free
definition through the spectral flow of Hessians, the conjugate
index is constructed using  an appropriate choice of coordinates
along the perturbed geodesic. The above theorem shows in
particular the independence of  conjugate index   from the various
choices involved in the construction.

   The classical Morse index theorem  is a  special case of Theorem
\ref{TEOREMA}. When  the metric is Riemannian $
\ispec(\gamma) \, = \,  \iMor(\gamma)$ while $\icon(\gamma ) =
\sum_{x \in ]0,1]} m(x).$

Our interest in  this problem was partially motivated   by the
questions about the stability of the focal index raised in
\cite{MPT}. However the main reason that lead us to the present
formulation of the Morse index theorem in terms of the spectral
flow is because of the relevance of this invariant to bifurcation
of critical points of strongly indefinite functionals found in
\cite{P}. Our  purpose is to combine the index theorem with
the results in \cite{P} in order to study bifurcation of perturbed
geodesics on semi-Riemanian manifolds. This will be done in
\cite{MPP}.

A few words about our definition of conjugate index are in order.
It is very close to other topological invariants that arise in
bifurcation theory. To some extent  it was suggested by the
approach to Hopf bifurcation in \cite{Iz} and  that of potential
operators in \cite{KK}  where  analogous  invariants are treated
in this way. The simplest topological invariants which detect
gauge anomalies are also of this form \cite{At}. Related ideas in
the context of Sturm-Liouville boundary value problems  can be
found in \cite{GST}.

Three beautiful papers \cite{ AtPS,GSig,RS1} have strongly
influenced the method of proof of Theorem \ref{TEOREMA}. This
proof has some interest in its own. It was found in trying to
understand the relation between the Morse index and regularized
determinants for families of boundary value problems discussed in
\cite{Le}. Although we did not quite succeed in this, yet we
believe that the proof  of our theorem  shed some light  on that
question.  Even in the case of Riemannian manifolds  it gives a
new proof of the classical Morse Index Theorem. Previous proofs
either used the variational characterization of eigenvalues of  a
self-adjoint operator (which cannot be used in the semi-Riemannian
case) or the homotopy properties of  Lagrangian Grassmannian
$\Lambda_n$. Here we substitute the later with the standard
properties of topological degree and operators in the trace class.

Other variants of our proof  can be  easily conceived.  A more
topological but less elementary one should be along the following
lines: the family of complex self-adjoint operators has an index
bundle in ${K}(S^2)\equiv { \Z}$ whose  first Chern class can be
easily related to the spectral index as in  \cite{BW}. On the
other hand one  can try to  relate this class to the conjugate
index by deforming the clutching function of this bundle to the
map $b_z$ in section $4.$ This will make the Morse index theorem
reminiscent of the Atiyah-Singer index theorem for  families of
{\em complex self-adjoint Fredholm operators\/}.

On the other hand,  after Floer's work,  results  of  type
spectral flow equals Maslov index became increasingly popular. See
for example \cite{FF},\cite{RS1} and the references there. The
paper of Nicolaescu \cite {N} places them into the realm of index
theorems for one parameter families of {\em real self-adjoint
Fredholm operators\/}. Index theorems for variational problems
with general Lagrangians and general boundary conditions, e.g.
focal points for geodesics starting from a sub-manifold  \cite{Z},
\cite{PPT}, can be also cast in the above form and indeed can
 be easily deduced as particular cases of  the main theorem in \cite{RS1}.

We consider here only  perturbations of the geodesic energy
functional by a potential because, as we mentioned above, this
  is the most general framework on which  the generalized Morse index
has an  intrinsic, geometric meaning. For more general
Lagrangians, as for example  in \cite{Z}, it still can be defined
but it depends on the choice coordinates. On the contrary our
theory   is general enough to cover trajectories of a mechanical
systems while keeping the geometric content of the geodesic case.

The paper is structured as follows: in section
\ref{sec:spectralflow} we shortly review the spectral flow of one
parameter families of quadratic functionals of Fredholm type. In
section \ref{sec:spectralindex} we give the variational
formulation of our geometric problem and  introduce the spectral
index. In section \ref{sec:conjugateindex}, using the orthonormal
parallel trivialization along the $p$-geodesic, we define the
conjugate index and compute it in terms of the associated Green
function. Theorem \ref{TEOREMA} is proved in section
\ref{sec:dimdelteorema}, while section \ref{sec:maslovindex} is
devoted to the relation with the Maslov index.

This is a second, slightly enlarged version, of the paper. The
first version circulated in the form of a  preprint from January
2003. We were unaware of the paper \cite{Z}  at that time.

Finally, let us point out that the main idea in \cite {PPT}( i.e., that
on semi Riemannian manifolds the spectral flow detects among
conjugate points those that are bifurcation  points), the one in
\cite{GPP} (that at degenerate conjugate points this can be
computed by means  of partial signatures) and the main idea in
this paper (that conjugate points can be also counted using
topological degree of a Wronskian) belong to the same project
and  were all conceived in discussions connected
with  the Ph.D. thesis of A. Portaluri at Polytechnic University of Turin.

\section{Spectral flow for paths of Fredholm quadratic
forms}\label{sec:spectralflow}

The informal description of the spectral flow of a path of
self-adjoint operators given in the introduction can be made
rigorous in many different ways. Beginning with \cite{AtPS},
several different approaches to this invariant appeared in the
literature. Here we will use the approach in \cite{P}. We will
need their construction in a slightly generalized form since our
goal is to give an intrinsic, i.e. coordinate free, construction
of the generalized  Morse index.

The tangent space to the manifold of paths on a semi-Riemannian
manifold has a natural Hilbertable structure but not a natural
Riemannian metric on it. Therefore we have to work directly with
paths of bounded quadratic (equivalently bilinear) forms arising
as Hessians of the energy  functional and not with self-adjoint
operators representing them with respect to a given scalar
product. Fortunately the spectral flow, being a topological
invariant, depends only on the path of quadratic forms and not, as
the name could misleadingly suggest, on the spectrum of the
operators representing the form. Below we will give the very
simple proof of this nontrivial fact. Then we will be able to
define spectral flow of a family of Fredholm  quadratic forms  on
a Hilbert bundle over the interval $ [0,1]$ which is the object
that intrinsically arises in  our framework.

Let $S,T$ be two invertible self-adjoint operators on a Hilbert
space $H$ such that  $S-T$ is compact. Then the difference between
spectral projections of $S$ and $T$ corresponding to a given
spectral set is also compact. Denoting with $E_{-}(\cdot)$ and
$E_{+}(\cdot)$ the negative and positive spectral subspace of an
operator,  it follows then that $E_{-}(S)\cap E_{+}(T)$ and
$E_{+}(S)\cap E_{-}(T)$ have finite dimension \cite{ P, ASS}.  The
{\em relative Morse index\/} of the pair $(S,T)$ is defined by
\[\irel(S,T)=\dim\left\{  E_{-}(S)\cap E_{+}(T)\right\}
-\dim\left\{ E_{+}(S)\cap E_{-}(T)\right\} .\]

It is easy to see that when the negative spectral subspaces of
both operators are finite dimensional $ \irel(S,T)$ is given by
the difference $ \iMor(S) - \iMor(T) $ between Morse indices.

A bounded self-adjoint operator $A$ is Fredholm  if $\ker A$ is
finite dimensional. The topological group $GL(H)$ of all
automorphisms of $H$ acts naturally on the space of all
self-adjoint Fredholm operators $\Phi_{S}(H)$  by {\em
cogredience\/} sending $A \in \Phi_{S}(H)$ to $S^*AS$. This
induces an action of paths in $GL(H)$ on paths in $\Phi_{S}(H)$.
It was shown in \cite[Therem 2.1]{P} that for any path $A \colon
[a,b]\to \Phi_{S}(H)$  there exist a path  $M
\colon[a,b]\rightarrow GL(H)$, and  a symmetry ${\mathcal J} (
{\mathcal J}^2= {\rm  Id})$ such that $M^*(t)A(t)M(t)= {\mathcal
J} +K(t) $ with $K(t)$ compact for each $ t\in [a,b].$

Let  $ A \colon [a,b]\to \Phi_{S}(H)$ be a path such that  $A(a)$
and $A(b)$ are invertible operators.
\begin{definition}\label{d3}
The  spectral flow of the path $A$  is the integer
\[\spfl(A,[a,b])\equiv \irel({\mathcal J}+K(a),{\mathcal
J}+K(b)),\] where $ {\mathcal J}+K$ is any  compact perturbation
of a symmetry cogredient with $A$.
\end{definition}

It follows from the properties of the relative Morse index that
the left hand side is independent of ${\mathcal J}, M$  and that
the above definition is nothing but a rigorous version of the
heuristic description of the spectral flow given in the
introduction \cite{P}.

The spectral flow $\spfl(A,[a,b])$ is additive and invariant under
homotopies with invertible end points. It is clearly preserved by
cogredience. For paths that are compact perturbations of a fixed
operator it  coincides with the relative Morse index of its end
points.

A Hilbertable structure on  a Hilbert space $H$ is the set of all
scalar products on $H$  equivalent to the given one.  A {\em
Fredholm quadratic form\/} is a  function $q\colon H\to \R$ such
that there exists a bounded symmetric bilinear form $b=
b_{q}\colon H\times H\to\RR$ with $q(u)=b(u,u)$  and with $ \ker b
$ of finite dimension. Here $\ker b = \{u : b(u,v) =0 \ \text{for
all}\  v \}$. The space $Q(H)$ of bounded quadratic  forms is a
Banach space with the norm  defined by $ \left\|  q\right\| =
\sup\left\{ \left| q(u) \right| : \left\| u\right\| =1\right\}$.
The set $Q_F(H)$ of all Fredholm quadratic forms is an open subset
of $Q(H)$ that is stable under perturbations by weakly continuous
quadratic forms. A quadratic form is called  {\em nondegenerate\/}
if the map $u\rightarrow  b_q (u,-) $ is an isomorphism between
$H$ and $H^*$. By Riesz representation theorem, for any choice of
scalar product $\left\langle \cdot,\cdot\right\rangle $ in the
Hilbertable structure,  $Q_F(H)$ is isometrically isomorphic to
$\Phi_{S}(H).$ Clearly this isometry sends the set of all
non-degenerate quadratic forms onto $GL(H)$. From the Fredholm
alternative applied to the representing operator it follows that a
Fredholm quadratic form $q$ is non-degenerate if and only if $\ker
b_q=0$.

A path of quadratic forms $q\colon [a,b]\rightarrow Q_{F}(H)$ with
nondegenerate end points  $q(a)$ and $q(b)$ will be called {\em
admissible\/}.

\begin{definition}\label{d4}

The spectral flow of an admissible path $q\colon [a,b]\rightarrow
Q_{F}(H)$  is given by
\[\spfl(q,[a,b])=\spfl(A_{q},[a,b])\] where $A_{q(t)}$ is the unique
self-adjoint operator such that $\left\langle
A_{q(t)}u,u\right\rangle ={q(t)}(u)$ for all $u\in H$.
\end{definition}

That this is independent from the choice of the scalar product in
a given structure follows from the invariance of the spectral flow
under cogredience. Indeed let $\left\langle
\cdot,\cdot\right\rangle_{1}$  be a scalar product equivalent to
$\left\langle \cdot,\cdot\right\rangle $ and let
$A_{q(t)},B_{q(t)}$ be such that $\left\langle
A_{q(t)}u,u\right\rangle =\left\langle B_{q(t)}u,u\right\rangle
_{1}=q(t)(u)\quad {\mbox{ for all }} u\in H$.

Denoting by $H_1$ the vector space $H $ endowed with the scalar
product $\left\langle \cdot,\cdot\right\rangle _{1}$, there exists
a positive self-adjoint operator   $D\colon H\to H_1$ such that $
\left\langle u,v\right\rangle _{1}=\left\langle Du,v\right\rangle
$ for all $u,v\in H.$ Therefore  $D ={\rm Id}^*$  where ${\rm Id}$
is considered as a map from $H_1$ into $H.$ Moreover we have that
$A_{q(t)}=DB_{q(t)}.$ Finally, by  invariance under cogredience we
get\[ \spfl(A_q,[a,b])=\spfl(DB_q,[a,b])=\spfl({\rm Id}^* B_q {\rm
Id} ,[a,b])=\spfl(B_q ,[a,b]),\] which is what we wanted to show.

We list below some properties of the spectral flow of an
admissible path of quadratic forms that we will use later. They
need not be proved here since they follow easily from the
representation formula in the Definition \ref{d4} and the
analogous properties of the spectral flow for paths of
self-adjoint  Fredholm operators  proven in \cite{P}.

\begin{itemize}
\item {\em Normalization\/}  Let $q \in C\left(  [a,b];
Q_{F}(H)\right)$ be such that $q(t)$ is non-degenerate for each
$t\in [a,b].$ Then $\spfl(q,[a,b]) =0$.

\item {\em Cogredience\/}  Let  $M\in
C\left([a,b];L(H_1,H)\right)$ be a path of  invertible operators
between the Hilbert spaces $H$ and $H_1$ and let $p$ be the path
of quadratic forms on $H_1$ defined by $p(t)(v) =
q(t)[M(t)^{-1}v]$. Then
$$\spfl(p,[a,b])=\spfl(q,[a,b]).$$

\item {\em Homotopy invariance\/} Let $h\in C\left(
[0,1]\times[a,b]; Q_{S}(H)\right)  $ be such that $h(s,t)$ is
non-degenerate for each $s\in [0,1]$ and $t=a,b.$ Then
$$\spfl(h(0,\cdot),[a,b])=\spfl(h(1,\cdot),[a,b]).$$

\item {\em Additivity.}  Let $c\in (a,b)$ be a parameter value at
which $q(c)$ is non-degenerate. Then
$$\spfl(q,[a,b])=\spfl(q,[a,c])+\spfl(q,[c,b]).$$
\end{itemize}

We will also need  a formula that leads to the calculation of the
spectral flow for paths with only   regular crossing points.

If a path   $q\colon [a,b]\rightarrow Q_F(H)$ is differentiable at
$t$   then the derivative $\dot{q}(t)$ is also a quadratic form.
We will say that a point $t$ is a crossing point  if $\ker
b_{q(t)}\neq \{0\} $, and we will say that the crossing point  $t$
is regular if the {\em crossing form\/} $\Gamma(q,t),$  defined as
the restriction  of the derivative  $\dot {q}(t)$ to the subspace
$\ker b_{q(t)}$, is nondegenerate. It is easy to see that regular
crossing points are isolated and that the property of having only
regular crossing forms is generic for paths in $Q_F(H)$. From
\cite[Theorem 4.1]{P} we obtain:
\begin{proposition}\label{crossform}
If all crossing points of the path are regular then they are
finite in number and
\begin{equation}\label{crformula1}
\spfl(q,[a,b])= \sum_i  \sgn \, \Gamma(q,t_i).
\end{equation}
\end{proposition}

A {\em generalized   family of Fredholm quadratic forms\/}
parameterized by an interval is a smooth function   $q\colon
{\mathcal H}\to \R ,$ where ${\mathcal H}$ is a Hilbert bundle
over $[a,b]$ and $q$ is such that its restriction $q_t$  to the
fiber ${\mathcal H}_t$  over $t$ is a Fredholm quadratic form. If
$q_a$ and $q_b$ are non  degenerate, we define the spectral flow
$\spfl(q) = \spfl (q,[a,b])$  of such a family $q$ by choosing a
trivialization $$M \colon [a,b] \times {\mathcal H}_a \to {\mathcal
H}$$ and defining
\begin{equation} \label{sflow2}
\spfl(q) = \spfl (\tilde{q},[a,b])
\end{equation}
where $\tilde{q}(t)u =q_t (M_tu).$

It follows from cogredience property that the right hand side of
\equ{sflow2} is independent of the choice of the trivialization.
Moreover all of the above properties hold true in this more
general case, including the calculation of the spectral flow
through a non degenerate crossing point given in Proposition
\ref{crossform} provided we substitute the usual derivative with
the intrinsic derivative of a bundle map.

\section{The spectral index}\label{sec:spectralindex}

Given a smooth $n$-dimensional manifold $M$, let $\Omega $ be the
manifold of all \mbox{$H^1$-paths} in $M$. Elements of $\Omega =
H^1 (I;M)$ are maps $\gamma \colon I \to M$ such that for any
coordinate chart $(U,\phi)$ on $M$ the composition $\phi \circ
\gamma \colon  \gamma^{-1} (U) \to \R^n$ belongs to $H^1
(\gamma^{-1} (U) ; \R^n)$. It is well known that $\Omega$ is a
smooth Hilbert manifold modelled by $H^1 (I;\R^n)$. We will denote
with $\tau\colon TM \to M $ the projection of the tangent bundle
of $M$ to its base, and by $H^1 (\gamma)$ the Hilbert space $H^1
(\gamma) = \{ \xi \in H^1 (I;TM) \, \colon  \, \tau\circ  \xi =
\gamma \}$ of all $H^1$-vector fields along $\gamma$. The tangent
space $T_{\gamma } \Omega$ at $\gamma$ can be identified in a
natural way with  $H^1 (\gamma).$ For all this the basic reference
is \cite{KL}.

By \cite[Proposition 2.1]{KL}, the map

\begin{equation}\label{fib}
\pi \colon  \Omega \to M\times M ; \, \quad \pi (\gamma ) =
(\gamma (0), \gamma(1))
\end{equation}
is a submersion and therefore for each $(p,q) \in M\times M$ the
fiber of $\pi$
\begin{equation}
\Omega_{p,q} = \{ \gamma \in \Omega \, :  \, \gamma (0)=p , \,
\gamma (1) = q \} \label{s2}
\end{equation}
is a submanifold of codimension $2n$ whose tangent space
$T_{\gamma } \Omega_{p,q} = \ker\,  T_{\gamma } \pi $ is
identified with  the subspace $H_0^1 (\gamma) $  of $H^1 (\gamma)
$ defined by
\begin{equation}
H^1_0 (\gamma ) = \{ \xi \in H^1 (\gamma ) \, : \, \xi (0)= \xi
(1) =0 \}. \label{s3}
\end{equation}
Since  $\pi$ is a submersion, it follows  that the family of
Hilbert spaces $H^1_0 (\gamma )$ form a Hilbert bundle  $TF (\pi)
= \ker \, T\pi $ over $\Omega$, called the  {\em bundle of
tangents along the fibers of\/} $\pi$.

To each pair $(g,V)$, where $g$ is a semi-Riemannian metric on $M$
and $V\colon I \times M\to \R$ is a smooth potential where
$I=[0,1]$, there is associated an {\em energy functional\/}
$E\colon \Omega \to \R$ defined by
\begin{equation}
E(\gamma ) = \int_0^1 \frac{1}{2} g( \gamma'(x) , \gamma'(x)) \,
dx - \int_0^1 V(x,\gamma (x)) \, dx. \label{s4}
\end{equation}
It is well known that $E$ is  a smooth  function and hence  so are
the restrictions $E_{p,q}$ of $E$ to $\Omega_{p,q}$. We will be
interested in the critical points of $E_{p,q}.$ The differential
of $E_{p,q}$ at a point $\gamma \in \Omega_{p,q}$ is given by the
restriction of $dE$ to $H^1_0 (\gamma ).$
      It is easy
to see that for  $\xi \in H^1_0 (\gamma)$
\begin{equation}
dE_{p,q} (\gamma ) [\xi ] =  \int_0^1 g( \Dds \xi(x),
\gamma'(x))dx - \int_0^1 g(\nabla V(x,\gamma (x)) , \xi (x) ) \,
ds . \label{s5}
\end{equation}
By standard regularity arguments one shows that if  $dE_{p,q}
(\gamma ) [\xi ] =0$ for all $\xi$ then $\gamma$  is  smooth and
then performing integration by parts one obtains that the
critical points of $E_{p,q}$ are precisely the smooth paths
$\gamma$  between $p$ and $q$ that verify the equation \equ{2.2}
of perturbed geodesics.

Let us recall that  if  $N $ is a Hilbert manifold and $n$ is a
critical point of a smooth function $f\colon N\to \R$ then the
{\it Hessian  of $f$}  at $n$  is the quadratic form  $h_n$ on
$T_n N $ given by $ h_{n} (v) = v( \chi (f))$, where $\chi $ is
any vector field defined on a neighborhood of $n$ such that
$\chi(n) =v.$  Through the identification of $ T_{\gamma}
\Omega_{p,q} $ with $H^1_0 (\gamma )$ a well know result in
Calculus of Variations  yields the Hessian of $E_{p,q}$ at
$\gamma.$  This is the quadratic form $h_{\gamma } \colon  H^1_0
(\gamma ) \to \RR$ whose associated bilinear form $H_{\gamma }
\colon  H^1_0 (\gamma)\times H^1_0 (\gamma ) \to \RR$  is given by
\begin{eqnarray} \label{s6}
  H_{\gamma } (\xi , \eta ) &=&\int_0^1 g( \Dds \xi (x) ,
  \Dds \eta (x)) \, dx \cr
\nonumber
  &-& \int_0^1 g( R(\gamma'(x) , \xi (x)
) \gamma'(x) + D_{\xi (x) } \nabla V(x,\gamma (x) ) , \eta (x) )
\, dx .
\end{eqnarray}

\begin{proposition}\label{fred} The form  $h_\gamma$ is a Fredholm
quadratic form. Moreover  $h_{\gamma}$ is non degenerate if and
only if $1$ is not a conjugate instant.
\end{proposition}

\proof We begin by constructing a Riemannian metric related to
$g$. Since $g$ is a non-degenerate symmetric form, we can split
$TM$ as direct sum of $T^+ M$ and $T^- M$ such that the
restriction of $g_\pm$ to $T^\pm M $ is positive definite and
negative definite respectively.

Let $j$ be the endomorphism of \mbox{$ TM =T^+ M \oplus T^- M$ }
given by $j(u^+ + u^- ) = u^+ - u^- .$ We define a new metric
$\bar g$ by $ {\bar g}(u,v)= g(ju,v) .$ Then $\bar g$ is a
Riemannian metric on $M$ and $j$ represents $g$ with respect to
$\bar g.$

The metric $\bar g$ induces a scalar product in $H^1_0(\gamma )$
given by
$$\langle \xi , \eta \rangle_{H^1_0} = \int_0^1
\bar{g}(\Dds \xi (x) , \Dds \eta (x)) \, dx .$$ By the very
definition of $ {\bar g} $ we have
\begin{equation}
\int_0^1 g(\Dds \xi (x) , \Dds \eta (x) ) \, dx = \langle{\mathcal
J_{\gamma}}\xi , \eta \rangle_{H^1_0} \label{s07}
\end{equation}
where ${\mathcal J_{\gamma}}(\xi)(x) := j(\gamma (x)) \xi (x)$
namely $\mathcal J_{\gamma}$ is pointwise $j$.

Clearly ${\mathcal J_{\gamma}}$ is bounded with ${\mathcal
J_{\gamma}}^2 =  I$, and hence the quadratic form
\begin{equation}
\label{defd} d_{\gamma } (\xi ) = \int_0^1 g(\Dds \xi (x) , \Dds
\xi(x)) \, dx
\end{equation}
  is non degenerate being represented by ${\mathcal
J_{\gamma}}\in GL(H^1_0(\gamma )).$ On the other hand $h_\gamma =
d_\gamma  - c_\gamma $  where
\begin{equation} \label{c}
c_\gamma (\xi) = \int_0^1 g\Big( R(\gamma'(x) , \xi (x) )
\gamma'(x) + D_{\xi (x) } \nabla V(x,\gamma (x) ) , \xi (x) \Big)
\, dx.
\end{equation}

The form  $c_\gamma$ is the restriction to $H^1_0 (\gamma )$ of a
quadratic form defined on the space $C_{\gamma}^0 (TM)$ of all
continuous vector fields over $\gamma.$ Since the inclusion $
H^1_0 (\gamma ) \hookrightarrow C_{\gamma}^0 (TM)$ is a compact
operator, it follows that $c_{\gamma }$ is weakly continuous and
therefore $h_{\gamma}$ is Fredholm being a weakly continuous
perturbation of a non degenerate form.

For the second assertion we notice that if $H_\gamma (\xi , \eta)
= 0$ for all $\eta \in H^1_0 (\gamma )$ then, again by regularity,
$\xi$ is smooth. Integrating by parts in \eqref{s6}, we
obtain that $\xi$  must verify the Jacobi equation \eqref{2.5}
with Dirichlet boundary conditions. The converse is clear.
Therefore $\ker \,h_\gamma = \{0\}$  if and only if  the instant
$1$ is not conjugate to $0$. \qedhere

\begin{rem} That the quadratic form $h_{\gamma}$ is Fredholm can be
proved without introducing the metric $\bar{g}$. However it is of
some interest to  notice that the use of $\bar{g}$ combined with a
parallel trivialization of the tangent bundle along $\gamma$
produces a concrete  construction of the abstract reduction of the
path of hessians to  a path of  compact perturbation of a symmetry
${\mathcal J}$ used in our definition of the spectral flow in
section 2.
\end{rem}

 From now on let $p, q \in M$ be fixed points and let $\gamma$ be a
$p$-geodesic from $p$ to $q$ such that $1$ is not a conjugate
instant.  Such a $p$-geodesic will be called {\em nondegenerate}.
In order to define the spectral index of a nondegenerate
$p$-geodesic $\gamma $ we will consider the  path induced  by
$\gamma$ on $\Omega.$

Namely, for each $t\in [0,1] $ let  $\gamma_t  \in \Omega$ be the
curve defined by $ \gamma_t (x) =\gamma (t\cdot x)$. Since
$\gamma$ is a critical point of $E_{p,q}$ it follows  from
\equ{2.2} that $\gamma_t$ is  a critical point of the functional
$E_t \colon \Omega_{p,\gamma (t)} \to \R$ defined by
\begin{equation}
E_t(\gamma ) = \int_0^1 \frac{1}{2} g( \gamma'(x) , \gamma'(x)) \,
dx - \int_0^1 t^2 V(tx,\gamma (x)) \, dx. \label{s44}
\end{equation}
In other words, each $\gamma_t$ is a $p$-geodesic for the
potential $V_t(x,m)= t^2V(tx,m).$

Let $h_t \colon  H^1_0 (\gamma_t) \to \R$ be the Hessian of $E_t$
at the critical point $\gamma_t$. By Proposition \ref{fred}, $h_t$
is degenerate  if and only if $1$ is  a conjugate instant for
$\gamma_t$. In particular  $h_1= h_\gamma$  is
non degenerate. Moreover $h_0$ is nondegenerate as well.  Indeed,
$\gamma_0 \equiv p$ is a constant path which is a critical point
of $E_0.$ An  $H^1$-vector field $\xi$ along $p$ is simply a path
$\xi \in H^1( I ; T_p(M))$ and hence $h_0(\xi)=d_{p} (\xi ) =
\int_0^1 g(\Dds \xi (x) , \Dds \xi(x)) \, dx $ which is
nondegenerate by the previous discussion.

Let us consider the canonical path of  $p$-geodesics
$\tilde\gamma \colon [0,1] \to \Omega$ defined by $\tilde\gamma
(t)  = \gamma_t .$  Clearly, the family of Hessians $h_t$, $0 \leq
t \leq 1$, defines a smooth function $h$ on the total space of the
Hilbert bundle $ {\mathcal H} =\tilde\gamma^{*} TF(\pi)$ over $[0
,1]$, that is a Fredholm quadratic form at each fiber and non
degenerate at $0$ and at $1.$ The spectral flow $\spfl(h)$ of such
a family  is well defined by \equ{sflow2} of the previous section.
\begin{definition} The {\em generalized Morse  index\/}
     $\ispec (\gamma )$ of  a
$p$-geodesic  $\gamma$ is the integer
\begin{equation}
\ispec(\gamma ) = -\spfl (h ) . \label{insp}
\end{equation}
\end{definition}
For pertubed geodesics on Riemannian manifolds the following
holds.
\begin{proposition}
If the metric $g$ is Riemannian  then  the Morse index
$\iMor(\gamma)$ (i.e. the dimension of the maximal negative
subspace of the Hessian of $h_{\gamma}$) is finite  and
\begin{equation}
\ispec (\gamma ) = \iMor(\gamma ) . \label{morse}
\end{equation}
\end{proposition}
          \proof  The first assertion is well known. It follows from
the fact that in the Riemannian case each $h_\gamma $ is a weakly
continuous perturbation of a positive definite form $d_\gamma.$
The dimension of the maximal negative subspace of this form
coincides with the dimension of the negative spectral space of any
self-adjont operator representing the form. But this subspace is
finite dimensional because the operator is  essentially positive,
i.e. compact perturbation of a positive one.

In order to prove  the second assertion we observe  that,  with
respect to the scalar product defined by formula \equ{s07} (with
${\mathcal J}\,=\,{\rm Id}$) on   $H^1_0 (\gamma_t),$ the form
$h_t$ is represented by an essentially positive operator of the
form ${\rm Id}-C_t$ with $C_t$ compact self-adjoint. By
\cite[Proposition 3.9]{P}  the spectral flow of a family of
essentially positive operators is the difference between the Morse
indices at the end points. Applying this to any trivialization of
$\mathcal H$ we have that
$$ \ispec (\gamma ) =  -\iMor(\gamma_0 )+\iMor(\gamma_1 ) =
\iMor(\gamma ).$$
  \qedhere

\section{The conjugate index}\label{sec:conjugateindex}

In this section we introduce a  topological invariant that counts
the  algebraic number  of conjugate points along a
 nondegenerate $p$-geodesic $\gamma$
 and which coincides with the expression
\equ{icon} in the case of a regular p-geodesic.

       Given a perturbed geodesic $\gamma$ we will use a particular
trivialization  of $\gamma^{*}(TM) $ by choosing  a $g-$frame
${\mathcal E}$ along $\gamma$ made by $n$  parallel  vector fields
$ \{ e^{1}, \ldots ,e^{n}\}.$ Here a $g-$frame means that   the
vector fields $e^ {i}$ are point-wise  g-orthogonal and moreover
$g(e^{i}(x),e^ {i}(x)) = \epsilon_i$, where $\epsilon_i =1$ for
$i=1, \ldots , n-\nu$, and $\epsilon_i =-1$ if $i\geq n-\nu+1$,
for all $x \in [0,1].$  Such a frame induces a trivialization
\begin{equation} M_{{\mathcal E}} \colon I\times \RR^n \to
\gamma^{ *} (TM) \label{triv}
\end{equation}
       of  $\gamma^{ *} (TM)$   defined by
$ M_{{\mathcal E}}(x, u_{1} \ldots u_{n}) =\sum_{i=1}^{n} u_{i} e^
{i}(x).$
\smallskip

Writing the vector field  $\xi$ along $\gamma$  as
         $\xi(x)=
\sum_{i=1}^{n} u_{i}(x)  e^ {i}(x),$
       inserting  the above expression in the equation \equ{2.5} of
Jacobi fields and taking $g$ product with $e^j$, we reduce the
Jacobi equation \equ{2.5}  to a linear second order system of
ordinary differential equations
$$ \epsilon_i u''_i (x) + \sum_{j=1}^{n} S_{ij}(x) u_{j}(x) =0 ; \quad
1\leq i\leq n,$$ where $S_{ij}= g(R(\gamma',e^i)\gamma' + D_{ e^i}
\nabla V(\cdot,\gamma) , e^j).$

Putting $ u(x)= (u_{1}(x), \ldots, u_{n}(x))$ and $S(x)=
(S_{ij}(x)) $, the above system becomes
\begin{equation}
           Ju''(x) + S(x)u(x)= 0
\label{jactriv}
\end{equation}
where $J$ is the symmetry
\begin{equation}\label{J}
J = \begin{pmatrix} {\rm Id}_{n-\nu} & 0 \cr 0 & -{\rm Id}_{\nu}

\end{pmatrix}
\end{equation}

Under the trivialization  $M_{{\mathcal E}}$,  the metric   $g$ on
$\gamma^{ *} (TM)$ goes into  indefinite product of index $\nu$ on
$\R^n$ given by $\langle u,v \rangle_\nu = \langle Ju,v\rangle$
where $\langle \cdot ,\cdot\rangle$ is the Euclidean scalar
product on $\RR^n.$ Since both $R(\gamma ' , -) \gamma'$  and the
Hessian $D_{-} \nabla V(\cdot,\gamma)$ are $ g$-symmetric
endomorphisms of $\gamma^{ *} (TM)$
%(see    \cite{O} Chapter 3, Lemma 49. )
it follows that the matrix $S(x)$ is symmetric.

Now let us apply the same argument  to each p-geodesic  $\gamma_t$
introduced in the previous section,  using the induced  parallel
g-frame \mbox{$ {\mathcal E}_t = \{ e^{1}_t, \ldots ,e^{n}_t\}$ }
with  \mbox{$ e^{i}_t(x)= e^{i}(t\cdot x).$} The corresponding
trivialization $M_{{\mathcal E}_t}$  transforms the equation for
Jacobi fields on $\gamma_t$ into  a one parameter family of second
order systems of ordinary differential equations
\begin{equation}
            Ju''(x) + S_t(x)u(x)= 0
\label{jactrivp}
\end{equation}
where $ S(t,x)\,=\,S_t(x) $ is smooth in $[0,1]\times[0,1]$ and
$S_t^*=S_t.$

By definition of $S_{ij}$ and $\gamma_t$ it follows  that
\begin{eqnarray*}
(S_t )_{ij}(x)&=& g(R(\gamma_t '(x),e_t^i(x))\gamma_t'(x) +
D_{e_t^i(x)} \nabla V(x,\gamma_t(x)),e_t^j(x)) =\cr
&=&t^2[g(R(\gamma'(t \cdot x),e^i (t\cdot x))\gamma(t\cdot x) +\cr
&+& D_{ e^i (t\cdot x)} \nabla V(t\cdot x,\gamma(t\cdot x)),
e^j(t\cdot x))] =  t^2 S_{ij} (t\cdot x)
\end{eqnarray*}
Hence we have
\begin{equation}\label{S}
S_t(x)=  t^2 S (t\cdot x).
\end{equation}
The trivialization $M_{{\mathcal E}_t}$ induces a one to one
correspondence  between  solutions $u= (u_{1}, \ldots, u_{n})$ of
the Dirichlet   problem
\begin{equation}
Ju''(x) + S_t(x)u(x)= 0 \quad u(0) =\,0\,= u(1) \label{bvp}
\end{equation}
and  Jacobi fields over $\gamma_t$  vanishing at $0,1.$ On the
other hand $\xi \to \xi_t$ is a bijection between  Jacobi fields
over $\gamma$ vanishing at $0$ and $t$ with the Jacobi fields over
$\gamma_t$ vanishing at $0,1.$  It follows  then  that $ t\in
(0,1] $ is a conjugate instant for $\gamma$ if and only if the
boundary value problem \equ{bvp} has  a nontrivial solution.
\smallskip

We now will take into account the complexified  problem \equ{bvp}
by considering  the operator $Ju'' + S_t(x)u$ acting on complex
valued vector functions $u\colon  I \to \C^n.$

Let ${\mathcal O}$ be the bounded  domain on the complex plane
defined by
\begin{equation}
\label{5.5} {\mathcal O} \, = \, \left \{ \, z = t+is \, \in \C \,
: \,0 < t <1, \, -1 < s  < \, 1 \right \}.
\end{equation}
For any $z = t+is \in\bar{ {\mathcal O}}$ let us  consider the
closed unbounded operator
$${\mathcal A}_z\colon  {\mathcal D}( {\mathcal A}_z)  \subset L^2
(I; \C^n ) \to L^2 (I ; \C^n )$$ with domain ${\mathcal D}(
{\mathcal A}_z) = H^2 (I;\C^n ) \cap H^1_0 (I; \C^n )$, defined by
\begin{equation}\label{A}
{\mathcal A}_{z} (u)(x) \, = \, Ju''(x) + S_z (x)u(x)
\end{equation}
where $ S_z (x) =S_t (x) + is \,{\rm Id}.$

The two parameter  family ${\mathcal A}_z$ of unbounded
self-adjoint  Fredholm operators of $H= L^2 (I ; \C^n )$ is a
perturbation of a fixed unbounded operator $Ju''$ by smooth family
$ {\mathcal S}\colon  \bar{{\mathcal O}} \to {\mathcal L}(H)$ of
bounded operators ${\mathcal S}_z$ defined by the pointwise
multiplication by the matrix  $S_z$ .

Let us consider now  the associated family of  first order
Hamiltonian  systems. Putting   $v= Ju',\ $ the equation  $
Ju''(x) + S_z u(x)=0$ becomes equivalent to
\begin{equation*}
\left\{\begin{array}{ll} u'(x) = Jv(x)  \\
v'(x) = - S_z (x) u(x)
\end{array}\right.
\qquad \forall x \in [0,1].
\end{equation*}

Taking $w=(u,v)\in \C^{2n}$, the above  system can be rewritten as
the complex Hamiltonian system
\begin{equation}\label{ham}
w'(x) = \sigma H_z(x) w(x)
\end{equation}
where \begin{equation} \label{sigma}\sigma =\begin{pmatrix} 0 &
-{\rm Id} \cr {\rm Id} & 0 \end{pmatrix}
\end{equation}
is the complex symplectic  matrix, while $H_z (x)$ is the matrix
defined by
\begin{equation}\label{acca}
H_z (x) = \begin{pmatrix} - S_z (x) & 0 \cr 0 & - J
\end{pmatrix}.
\end{equation}
Let $\Psi_z(x)$ be the fundamental solution of \equ{ham}. The
matrix $\Psi_z(x)$ is the unique solution of the Cauchy problem
\begin{equation}\label{5.7}
\left\{\begin{array}{ll} \Psi_{z}'(x) = \sigma  H_{z}(x) \Psi_{z}
(x)  \quad
x \in [0,1]\\
\Psi_{z} (0) = {\rm Id}.
\end{array}\right.
\end{equation}
Consider the block decomposition of $\Psi_z(x)$
\[
\Psi_z (x)\,= \,  \begin{pmatrix} a_z (x) & b_z(x) \cr c_z(x) &
d_z(x)\end{pmatrix}\] and let $b_z = b_{z}(1)$ be the upper right
entry in the block decomposition of $\Psi_z(1)$.

The matrix  $b_z$ can be also described directly  in terms of the
solutions of \equ{jactrivp} as follows: for each $i, \, 1\leq
i\leq n,$ let  $u_i$ be a solution of the initial value problem
$Ju''(x) + S_z u(x)=0, \ \  u_i (0)=0, u^{\prime}_i(0) = e_i, $
where $ e_1,\dots , e_n $ is the canonical basis of $\C^n$ then
$b_z = ( u_1(1), \dots, u_n(1))^T$.

Our definition of the conjugate index is based on the following
elementary  observation (compare \cite[Lemma 1.5]{M})

\begin{lemma}\label{5.8}
The following three statements are equivalent:
\begin{enumerate}
\item[(i)] $\ker{\mathcal A}_{z}\neq \{0\}$; \item[(ii)] $ Im (z)
=0 $ and $t = Re (z)$ is a conjugate instant; \item[(iii)] $ \det
\, b_{z}\, = \, 0$.
\end{enumerate}
\end{lemma}

\proof  Since the spectrum of  self-adjoint operator is real,
$\ker{\mathcal A}_{z}\neq 0$ can occur only at $z= t+is$ with
$s=0$. But functions belonging to the kernel of  ${\mathcal
A}_{t}$ are precisely the solutions of the boundary value problem
\equ{bvp}. Therefore $t$ must be a conjugate instant. The converse
is clear. Hence the equivalence between $(i)$ and $(ii)$ is
proved.

Now let $u\in \ker{\mathcal A}_z $. From the block decomposition
of $ w(1)=\Psi_z (1)w(0)$ and the boundary conditions
$u(0)=0=u(1)$ we get  $b_zu'(0)\,=\,0.$ Thus, if $\det \,b_{ z}
\not= 0$ we have $u'(0) = 0$ and hence $u \, \equiv \, 0$, which
yields $\ker \, {\mathcal A}_{z } \, = \, \{0\}$. On the other
hand, if $\det\, b_{z} \,=\,0$, taking  $0\not= v_0 \, \in \ker\,
b_{z}, \ \ w_0 = ( 0, v_0)$ and $u (x) $ equal to the first
component of $\Psi_z (x) w_0,$ we have $ u(x) \not\equiv 0$ and $
u \in \ker \, {\mathcal A}_{z }$. This proves the equivalence
between $(ii)$ and $(iii)$.\qedhere

Consider now   $\rho  \colon  \overline {\mathcal O} \rightarrow
\, {\C}$ given by   by $\rho(z) \,=\, {\rm det}\,b_{z}$. Since
       $\ker \mathcal A_{0 }=\ker \mathcal A_{1 }=\{0\} $
  it follows that $0 \, \not\in \, \rho \,(\partial \, {\mathcal O})$.
Under this conditions the Brouwer degree $\deg \, (\rho \, , \,
{\mathcal O} \, , 0)$ of the map $\rho$ in ${\mathcal O}$ with
respect to $0$ is defined \cite{Deim}. Brouwer's  degree is a
topological invariant that counts with multiplicities the number
of zeroes of $\rho$ in ${\mathcal O}$ and since the zeroes of
$\rho$ correspond to conjugate instants of the $p$-geodesic we
make the following:

\begin{definition}\label{5.10}
The conjugate index $\icon\,(\gamma)$ of a nondegenerate
$p$-geodesic $\gamma $ is the integer
\begin{equation}\label{mu}
\icon \, ( \gamma)\,= -\deg \, ( \rho , \, {\mathcal O} , \, 0).
\end{equation}
\end{definition}

That the right hand side is independent from the choice of
${\mathcal O}$ follows from the excision property of degree.

The algebraic multiplicity of an isolated, but not necessarily
regular, conjugate instant $t_0$ can now be defined by
\begin{equation}\label{multip}
\icon(\gamma, t_0)=\,-\deg \, ( \rho , \, {\mathcal O}' , \, 0),
\end{equation}
where ${\mathcal O}'$ is any open neighborhood of the point $
z_0\,=\,(t_0,0)$ in ${\mathcal O}$ containing no other zeroes of
the map $\rho.$

If all conjugate instants are isolated then $\icon \, ( \gamma)\,
$ is the sum of the algebraic multiplicities of the conjugate
instants.

In section \ref{sec:maslovindex} we will show that if $t_0$ is  a
regular conjugate instant, then $  \icon\,(\gamma,t_0) $ coincides
with $ \sgn\,\big( g|_{{\mathcal I}[t_0]^{\bot}}\big)$.

If all data including  are analytic, the  multiplicities  of the
conjugate points can  be computed using an algorithm for
computation of degree of a polynomial plane vector vector-field
whose origins can be traced back to Kroenecker and Ostrogradski.
Below we shortly sketch a calculation of this type in a special
case.

Assuming $M, g $ and $V$ analytic, by uniqueness  of solutions of
ordinary differential equations we have that the restriction of
$\rho$ to the real axis  is a real analytic function which does
not vanish identically and hence has only isolated zeroes.
Therefore $\rho$ has a finite number of zeroes of the form $z_i
=(t_i,0)$ with  $t_1 < t_2 < \ldots < t_k$. Taking isolating
neighborhoods ${\mathcal O}_i$  as before we get
\[ \icon\,(\gamma)= \sum_{i=1}^k\icon\,(\gamma, t_i)
   =  -\sum_{i=1}^k \deg \, ( \rho ,\, {\mathcal O}_i, \, 0).\]
Fix a conjugate instant $t_j$  and let $ z_j =(t_j,0).$ With our
assumptions  $\rho $ is a  real analytic map from $\C \cong
\R^{2}$ into itself. Let $P$ and $Q$  be the non vanishing
homogeneous polynomials of lowest degree (respectively $m$ and
$n$) that arise in the  (real) Taylor series of the map $\rho$ at
the point $z_j.$

Then we have
\begin{equation}
\label{pol} \rho (z)= \rho (t,s) = (P (t-t_j ,s)  + f(t,s), Q
(t-t_j , s) + g(t,s))
\end{equation}
\noindent  with $| f(z)| = O(|z|^{m+1}) $ and $| g(z)| =
O(|z|^{n+1})$  for  $z$ ranging on a bounded set.

We will present the computation of  $\icon( \gamma, t_j)$ upon the
extra assumption:
\begin{enumerate}
\item[$(H_1)$] The point  $z_j$ is  an isolated zero of the
homogeneous map $\eta= (P,Q)$.
\end{enumerate}
The general algorithm for the degenerate case can be found in
\cite[ Chap. I ß15]{KZ} and \cite[Appendix]{KPPZ}.

       We first show
that the local multiplicity of $\rho $ and $\eta$ at $z_j$ are the
same using a well known argument. Without loss of generality we
can assume $z_j=0$. By homogeneity of $P$ and $Q$ and compactness
of the unit circle we can find two positive constants  $A$ and $B
$ such that  for each $z\neq 0$  either $  A |z|^m  < P(z) $ and
$| f(z)| < B(|z|^{m+1})$ or $  A |z|^n  < Q(z) $ and  $| g(z)| <
B(|z|^{n+1}).$ It follows from this inequalities that the homotopy
   $h(\lambda, z) = (P (z)  + \lambda f(z), Q (z ) + \lambda g(z)) ; \,
0\leq \lambda \leq 1$ does not vanishes on the boundary of the
circle with center at $0$ and radius $  R= A/B$. Therefore $ \icon
( \gamma, t_j) = - \deg \, ( \rho ,\, {\mathcal O}_i, \, 0)= -\deg
\, ( \eta ,\, B(0, R) ,\, 0)$. The right hand side can be computed
directly from the coefficients of $P$ and $Q$ as follows: we first
observe  that $0$ is an isolated zero for $\eta $ if and only if
$N_0 ( s ) = P (1,s)$ and $N_1 (s) = Q (1,s)$ do not have common
real root and $ |P (0,1) | + |Q (0,1) | >0$.

Assuming $m\geq n,$ let us consider polynomials $N_0 $,
$N_1,\ldots , N_l $, defined inductively by   the Euclidean
algorithm, namely $N_{i+1} (s) $ is the rest of the division of
$N_{i-1} (s) $  by $N_{i} (s)$. The final polynomial $N_l (r)$ is
the greatest common divisor of the polynomials $N_0$ and $N_1$.

Choose  $r$ in $\R$ such that $r$ is not root of any of the
polynomials $N_{i}$. We denote by $m(r)$ the number of sign
changes of the corresponding values $N_i (r).$ The number $m(r)$
becomes constant for $r$ sufficiently large. We denote with
$m_{+}$ this  value. Analogously, we denote by $m_{-}$ the common
value of $m(r)$ for $r$ negative and of sufficiently large
absolute value. From the above discussion and \cite[Theorem
10.2]{KZ} we have the following formula for the conjugate index at
a degenerate point.

\begin{lemma}\label{kr-zab}

If  $M, V$ are analytic and if $(H_1)$ holds at  a conjugate point
$t_j$ of $\gamma,$ then
$$\icon ( \gamma, t_j) = -\left( 1+ (-1)^{m + n} \right) \frac{m_{+} -
m_{-}} {2} .$$
\end{lemma}
\smallskip

In \cite{MPT} the authors defined the conjugate index using  the
formula \equ{icon} irrespectively whether the geodesic is regular
or not. They constructed  an example of a geodesic with an
isolated conjugate point and such that the associated Maslov index
does not coincide with the equation \equ{icon}. The correct
multiplicity at such a point is not given by  $\sgn\,
g|_{{\mathcal I}[t]^{\bot}}$. It can be computed  either using the
algorithm for calculation of the degree of an analytic vector
field in our approach  or  a well known formula for the Maslov
index  in terms of  partial signatures  in the approach chosen in
\cite{GPP}.  Let us remark in this respect that the multiplicity
of a conjugate point, as defined in \equ{multip}, coincides with
the analogous multiplicity defined in terms of the partial
signatures in \cite{GPP}. This follows immediately from
Proposition \ref{mas=spec}. However we were unable to find any
direct relation between the invariants $m_{\pm} $ arising in the
above lemma with the partial signatures at the given point.

We close this section with a proposition which provides a way to
compute the number defined by the formula \equ{mu} in terms of the
trace of the Green kernel of \equ{A}. This is  essentially  a
known result \cite{Le,GGK}. We include the proof here for the sake
of completeness.

Let us recall first that a compact operator  $K$  is said to be of
{\em trace class\/} if the series of the square roots of
eigenvalues of $K^{*}K$ is  convergent. The trace class ${\mathcal
T}$ is a bilateral ideal contained in the ideal of all compact
operators ${\mathcal K}.$ There is a well defined linear
functional ${\rm Tr}$ on ${\mathcal T}$ which has the usual
properties  of the trace. In particular if both $ AB$ and $BA$ are
in ${\mathcal T}$ then ${\rm Tr} AB = {\rm Tr} BA.$

Our calculations below can be better formalized using
operator-valued differential one-forms.  By a slight abuse of
notation we will denote by  ${ d\mathcal A}_z$ the differential of
the bounded part ${\mathcal S}_z$ of the family ${\mathcal A}_z$
, and consequently we will denote by  $d{\mathcal A}_{z}\,
{\mathcal A}^{-1}_{z} $ the operator valued  one-form given by
\begin{equation}\label{tr}
d{\mathcal A}_{z}\, {\mathcal A}^{-1}_{z} =
\partial_t {\mathcal S}_{z}{\mathcal A}_{z}^{-1}\, dt +
i \partial_s {\mathcal  S}_{z}{\mathcal A}_{z}^{-1}\, ds =
\partial_t{\mathcal S}_t {\mathcal A}_{z}^{-1}\, dt +  i {\mathcal
A}_{z}^{-1}\, ds
\end{equation}
defined on the set $ \{z\in \bar{\mathcal O}: {\mathcal A}_{z}
{\rm \,has\, a \, bounded\,  inverse} \}$.

This notation incorporates  the action on the left (resp. right)
of an operator valued function on a  one-form  in the natural way,
by multiplying on the left (or right) the coefficients of the
form. In the same vein, the trace  ${\rm Tr}\,\theta $  of an
operator valued one form $\theta =E\,dt +F\,ds$ is the complex
valued one form ${\rm Tr}\,E\,dt + {\rm Tr}\,F\,ds.$

\begin{proposition}\label{prop5.1}
          The form   $d{\mathcal A}_{z}\, {\mathcal A}^{-1}_{z} $ is a
trace class valued form and
\begin{equation}\label{tr1}
\icon(\gamma) \, = -\, \frac{1}{2 \pi i}\, \int_{\partial
{\mathcal O}} {\rm Tr}\,d{\mathcal A}_{z}\, {\mathcal A}^{-1}_{z}
.
\end{equation}
\end{proposition}

\proof  By Lemma \ref{5.8}, ${\mathcal A}_z$ has a bounded inverse
${\mathcal A}_z^{-1}$ for all $z \in \partial {\mathcal O}.$
Moreover, the operator ${\mathcal A}_z^{-1}$ is an integral
operator of the form
\begin{equation}
\label{ir} {\mathcal A}_z^ {-1}(u)(x)\, =\, \int_{0}^{1} K_z (x\,
,\,y)\, u(y)\, dy
\end{equation}
with the Green-kernel $K_z (x \, , \, y)$ given by
\begin{equation}\label{5.2}
K_z (x \, , \, y  ) = C \widetilde K_z (x\, , \, y ) D^*,
\end{equation}
where $C \, = \, ( I \, , \, 0),$ $D \, = \, (0 \, , \, I)$ and
$\widetilde{K}_z (x,y)$ is the $2n \times 2n$ matrix defined by
\begin{equation}
\widetilde{K}_{z}(x,y)=\left \{\begin{array}{ll}  - \Psi_z(x) P_z
\Psi_z(y)^{-1} & \mbox{ $0  \leq  x < y  \leq  1$} \cr
   \Psi_z(x) \,(I\, - \, P_z)
\Psi_z(y)^{-1}  \,
             & \mbox{$
0  \leq  y < x  \leq  1$}
\end{array}\right.
\end{equation}
with
\begin{equation}
P_z \, =\, \begin{pmatrix} 0 &
           0 \cr b_z^{-1} & 0
           \end{pmatrix}\, \Psi_z (1)
\end{equation}
(see \cite[Chapter XIV, Theorem 3.1]{GGK}).

Since  the kernel in  \equ{5.2}  is of class $C^{0,1},$ it follows
from a well known theorem of Fredholm that  the operator $
{\mathcal A}_z^{-1}$ is of trace class and therefore  form
$d{\mathcal A}_z {\mathcal A}_z^{-1}$ is
      ${\mathcal T}-$valued.

In order to prove the formula \equ{tr1} we will first calculate  $
{\rm Tr} \,  \,d{\mathcal A}_{z} \, {\mathcal A}_{z}^{-1}\,$ using
the fact that  the trace of an integral operator belonging to the
trace class can be computed integrating the trace of its kernel
\cite{GGK}. For $z \in\partial {\mathcal O}$ from  \equ{sigma} and
\equ{acca} we get
\[{\rm Tr} \,  \,d{\mathcal A}_{z} \, {\mathcal A}_{z}^{-1}\, =
\int_0^1 {\rm Tr} [ dS_z (x)\, K_{z}(x,x) ]dx=- \int_0^1 {\rm Tr}
[\sigma  dH_z(x) \widetilde K_z(x,x)]dx.\] On the other hand,
\begin{eqnarray*}
&-& {\rm Tr}  [ \,\sigma \,d{H}_{z}(x) \, \, \tilde K_{z}(x \, ,
\, x)\,] =\,{\rm Tr} \, [ \, \sigma d{H}_{z}(x) \, \Psi_{z}(x) \,
P_z\, \Psi_{z}^{-1}(x)]=\cr &=& {\rm Tr} \, [ \, d{\Psi}'_{z}(x)\,
P_z\,\Psi_{z}^{-1}(x) \, - \, \sigma \, H_z (x) \, d \Psi_z (x) \,
P_z \, \Psi_z^{-1} (x)  \, ]= \cr &=&{\rm Tr} \, [
\,d{\Psi}'_{z}(x)\, P_z\,\Psi_{z}^{-1}(x)\,-\,\Psi'_{z}(x) \,
\Psi_{z}^{-1}(x)\, d{\Psi}_{z}(x)\,P_z\, \Psi_{z}^{-1}(x)\,]= \cr
&=&\frac{d}{dx} \,{\rm Tr} \, [ \, d{\Psi}_{z}(x)\,
P_z\,\Psi_{z}^{-1}(x)\,].
\end{eqnarray*}
 From this, integrating  in $x$,  a direct computation yields
$${\rm Tr} \,  d {\mathcal A}_z {\mathcal A}_z^{-1} =
      {\rm Tr} \,  [ \,d{\Psi}_{z}(1)\,
P_z\,\Psi_{z}^{-1}(1)\,]=
            {\rm Tr}\,  d b_z  b_z^{-1} = \, d \, log
\rho (z),$$ since $ {\rm Tr}\,  d b_z  b_z^{-1}$ coincides with
the logarithmic differential of ${\rm det} \, b_z.$
\smallskip

Integrating over $\partial {\mathcal O}$ we finally obtain
$$\frac{1}{2 \pi i}\, \int_{\partial {\mathcal O}} \,{\rm Tr} \,
d{\mathcal A}_{z}\, {\mathcal A}_{z}^{-1}\, =\,\frac{1}{2 \pi
i}\,\int_{\partial {\mathcal O}} \,  d\, log \,  \rho (z) \, ,$$
    which is precisely the degree $\,\deg
\, ( \rho , \, {\mathcal O} , \, 0)$ by the well known formula
relating the degree of a map on the open set ${\mathcal O}$ with
the winding number of its restriction to  the boundary (see
\cite[Chapter 1, Section 6.6]{Deim}).

\section{Proof of Theorem \ref{TEOREMA}}\label{sec:dimdelteorema}

Let  $ \tilde\gamma  \colon  [0 ,1] \to \Omega; \ \
\tilde\gamma(t) =\gamma_t$ be  the canonical path defined in
section \ref{sec:spectralindex} and let $h\colon {\mathcal H} \to
\R$ be the generalized family of quadratic forms whose restriction
$h_t$ to the fiber $H^1_0 (\gamma_t)$ of the Hilbert bundle
${\mathcal H}= \tilde\gamma^* TF(\pi) $ is the Hessian of $E_t$ at
$\gamma_t$. The parallel trivialization $M_{{\mathcal E}_t}$ of
$\gamma_t^* TM$ defined by formula \equ{triv} induces a
trivialization of ${\mathcal H},$ under which $h$ is transformed
into the family of Fredholm quadratic forms $\hat h_t$ on  $H^1_0
([0,1] ; \R^n)$ given by $\hat h_t (u) = h_t (\sum_{i=1}^{n} u_{i}
e_t^ {i})$. Using the computation of the previous section we
obtain
\begin{equation}\label{hess}
\hat h_t (u) =  \int_0^1 \langle J u'(x) , u'(x) \, \rangle dx -
\int_0^1 \langle S_t(x) u(x) , u(x)\rangle \, dx
\end{equation} where $J$ is given by \eqref{J} and $S_t$
is the smooth family of symmetric matrices introduced in \equ{S}.

By definition of the generalized Morse  index of $\gamma$, we have
\begin{equation}
\isf(\gamma) =- \spfl (h) = -\spfl (\hat{h} , [0 , 1]).
\label{eq:t5}
\end{equation}
We will reduce the calculation of $\spfl (\hat{h},[0,1])$ to that
of a path having only regular crossing. In order to obtain this we
will apply a perturbation result of Robbin and Salamon in
\cite{RS1} to the path of operators $ \tilde{{\mathcal
A}}=\{\tilde {{\mathcal A}}_t \}_{\{t \in [0 , 1]\}} $ where \-
$\tilde {{\mathcal A}}_t\colon {\mathcal D}( \tilde {{\mathcal
A}}_t) \subset L^2 (I ; \RR^n ) \to L^2 (I ; \RR^n )$ is the
closed, real self-adjoint operator defined on ${\mathcal D}(
\tilde {{\mathcal A}}_t) = H^2 (I;\RR^n ) \cap H^1_0 (I; \RR^n )$
by
\begin{equation}\label{AR}
\tilde {{\mathcal A}}_{t} (u)(x) \, = \, Ju''(x) + S_t (x)u(x).
\end{equation}

Notice that the restriction $\{{\mathcal A}_t;  t \in [0 , 1]\}$
of the previously defined  family $\{{\mathcal A}_z\}$ given by
formula \eqref{A} to the real axis is nothing but the
complexification of the path $\tilde {{\mathcal A}}.$

Since  $\tilde {{\mathcal A}}_t = JD^2 +{\mathcal S}_t $ is  a
compact differentiable perturbation of $ JD^2$, it verifies all
the assumptions in \cite[Theorem 4.22]{RS1} and hence there exist
a $\delta >0$ arbitrarily close to zero such the path
$\tilde{{\mathcal A}}_t^\delta = \tilde {{\mathcal A}}_t +\delta I
$ has only regular crossing points. Regular crossing for paths of
unbounded operators has the same  meaning  as
  for paths of quadratic forms. Namely, $t_0 $ is a regular
crossing point if the crossing form $\Gamma(\tilde{{\mathcal
A}}^\delta, t_0),$ defined as the restriction of the quadratic
form
\begin{equation} \label{eq:crossjac1}
\langle \dot {\mathcal S}_t u,u\rangle_{L^2} = \int_0 ^1 \langle
\dot S_t (x)u(x),u(x)\rangle\, dx
\end{equation}
to $\ker  \tilde {{\mathcal A}}^\delta_{t_0}$, is non degenerate.
Regular crossing points  are isolated  and thus $ \Ker \tilde
{{\mathcal A}}_{t}^\delta \neq\{0\}$ only at a finite number of
points $0 < t_1 < \ldots < t_k <1.$

Let   $\hat h_t^{\delta} (u)  = \hat h_t (u) +\frac 12 \delta
\left|\left| u\right|\right| _{L^2}$. From equation \equ{hess},
using once more integration by parts, we obtain that $ \Ker \,
\hat h_t^{\delta} = \Ker \, \tilde {{\mathcal A}}_{t}^\delta$.
Moreover, by the very definition of the crossing form, at each
crossing point $t_j$  the forms $\Gamma(\hat h ^{\delta}, t_j) $
and $\Gamma(\tilde{{\mathcal A}}^\delta, t_j)$ coincide. Taking
$\delta$ small enough, using the homotopy invariance of spectral
flow and Proposition \ref{crossform} we obtain
$$ \isf(\gamma ) = -\spfl ({\hat h}^\delta , [0 ,1] ) =
-\sum_{j=1}^k \sgn \, \Gamma({\hat h}^\delta , t_j ) =
-\sum_{j=1}^k \sgn \, \Gamma (\tilde {{\mathcal A}}^\delta , t_j
).$$

On the other hand, the complexified path $\{{\mathcal
A}^{\delta}_t ; \  t \in [0 , 1]\} $ has the same  crossing points
as $ \tilde {{\mathcal A}}^\delta$ because  $\Ker {\mathcal
A}_{t}^\delta$ is the complexification  of $ \Ker \tilde
{{\mathcal A}}_{t}^\delta$ and $ \sgn\, \Gamma({\mathcal A}^\delta
, t_j )\,= \, \sgn\, \Gamma(\tilde {{\mathcal A}}^\delta , t_j )$
because a real symmetric  matrix and its complexification  have
the same eigenvalues. This together with the previous calculation
gives

\begin{equation}
\isf(\gamma )=  -\sum_{j=1}^k \sgn \, \Gamma ({\mathcal A}^\delta
, t_j ).
   \label{eq:t06}
\end{equation}

Now let us compute the perturbation on the other side of the
equality in Theorem \ref{TEOREMA}. Let $ b^\delta_{z}(x)$ be the
upper right entry in the block decomposition of the fundamental
matrix $ \Psi^\delta_{z}(x)$ associated to the perturbed family $
{\mathcal A}^\delta_{z }$  and let $ \rho^\delta (z) \, = \, {\rm
det}\, b^\delta_{z}(1)$. By the continuous dependence of the flow
with respect to the data, we can take $\delta$ small enough so
that $ | \rho^\delta (z) -\rho(z)| < \inf _{z \in
\partial {\mathcal O}}|\rho(z)|.$
Then  $h(t,z) \, = \, t \rho^\delta (z) +(1-t) \rho(z)$ is an
admissible homotopy and therefore by the homotopy invariance of
degree
\begin{equation} - \icon(\gamma) \, = \,\deg ( \rho ,
{\mathcal O} , 0 ) \, = \, \deg(  \rho^{\delta} , {\mathcal O} ,
0). \label{eq:dim}
\end{equation}
Taking closed disjoint neighborhoods ${\mathcal D}_j$ of $(t_j,0)$
in $ {\mathcal O}$ with  piecewise smooth boundary from the
additivity of degree and from Lemma \ref{5.8} we obtain
\begin{equation}
-\icon(\gamma) \, = \,\sum_{j=1}^{k}\deg (  \rho^\delta ,
{\mathcal D}_j , 0 ) = \sum_{j=1}^{k} \,\frac{1}{2 \pi i}\,
\int_{\partial {\mathcal D}_j}  {\rm Tr}\,d{\mathcal A}^\delta_{z}
( {\mathcal A}^\delta_z)^{-1} . \label{eq:t6}
\end{equation}
Comparing  \eqref{eq:t06} with \eqref{eq:t6}, we see that the
proof of the Theorem \ref{TEOREMA} will be complete if we prove
the identity
\begin{equation}\label{1}
\frac{1}{2 \pi i}\, \int_{\partial {\mathcal D}_j} {\rm Tr}\,
d{\mathcal A}^{\delta}_{z}\,({\mathcal A}^{\delta}_{z})^{-1}\, =
\, \sgn \,\Gamma ({\mathcal A}^\delta , t_j ) \quad {\text{for any
}}\  j=1, \ldots , k.
\end{equation}

The rest of the section will be devoted to the proof of \eqref{1}.
\begin{convention}
On the basis of the previous discussion we  can assume that
${\mathcal A}_{t }$ has only  regular crossing points and drop
$\delta $ from our notations.
\end{convention}
The idea of the proof is as follows:  in a standard way one can
find a path $P_t$ of finite rank projectors defined on a
neighborhood of $t_j$  such that  $P_t$ reduces  ${\mathcal A}_{t
}$ and hence also ${\mathcal A}_{z }$ with  $ \Re z =t.$ By a
well-known theorem in \cite{K},  after a $t$-dependent unitary
change of coordinates, the family can be locally reduced by a
single projector. Using this reduction equation \eqref{1} follows
from the corresponding statement in finite dimensions where it can
be shown to be true by an  elementary calculation. This  works
well for bounded operators.  However in our setting a problem
arises.  Under change of coordinates the transformed family is not
any more a bounded perturbation of a fixed operator and there is a
problem with the definition of the one form in \equ{1}. We will
avoid all  the technicalities related to the differentiability of
general families of unbounded operators by rewriting all the forms
in terms of bounded operators.  For this we will deal
simultaneously with ${\mathcal A}_{z }$ as closed operators and
also as bounded operators with respect to the graph norm on the
domain. Notational  ambiguities usually arise in such a situation.
Hence, even at cost of being clumsy, we will carefully distinguish
both operators on our notation and we will do the same with the
corresponding families of spectral projectors.

We will denote  by $ W$ the Hilbert space $ H^2 (I ; \C^n) \cap
H_0^1 (I ; \C^n ) $  endowed with the graph norm of $Ju''$ and, as
before, $H$ will denote $L^2 (I; \C^n ).$ Let ${\widetilde j}$ be
the inclusion  of $ W$ into $H$.  The family $ A_{z }\, = \,
{\mathcal A}_{z } \circ {\widetilde j}$ is a family of {\em
bounded\/} operators. Moreover, since the operator ${\widetilde
j}$ is compact we have that  ${\mathcal A}_z^{-1} ={\widetilde j}
\circ A_z^{-1}$ is a compact operator whenever $ A_z^{-1}$ exists
and is bounded. On the other hand, being $dA_z = d {\mathcal A}_z
\circ {\widetilde j}$, we have also that
\begin{equation}
d{\mathcal A}_z \circ {\mathcal A}_z^{-1} = d A_z \circ A_z^{-1} .
\label{eq:Aidjac1}
\end{equation}
For a fixed $j$,  choose a positive number $\mu> 0$ such that the
only point in the spectrum of ${\mathcal A}_{t_j}$ in the interval
$[-\mu , \mu ]$ is 0 and then  choose $\eta$ small enough such
that neither $\mu$ nor $- \mu$ lies in the spectrum of $ {\mathcal
A}_{t }$ for $ |t-t_j|<\eta.$ For such a $t$, let $P_{t}$ be the
orthogonal projection in $H$ onto the spectral subspace associated
to the part of the spectrum of ${\mathcal A}_{t}$ lying in the
interval $[-\mu , \mu ].$ Then ${\mathcal A}_{t}P_t = P_t
{\mathcal A}_{t}$ on the domain of  ${\mathcal A}_{t}.$ In other
words $P_t$ reduces the operator ${\mathcal A}_{t}$.

 From  the  integral representation of $P_t$  given by
\begin{equation}\label{eq:proiettorijac1}
P_{t} \, = \, \frac{1}{2 \pi i} \int_{{\mathcal C}} \, ( {\mathcal
A}_{t} \, + \, \lambda \, {\rm Id} )^{-1} \, d \lambda,
\end{equation}
where ${\mathcal C}$ is a symmetric  curve in the complex plane
surrounding the spectrum in $( - \mu  , \mu ),$ one can show that
the projector $P_t$ factors through ${\widetilde j}.$ Indeed,
defining ${R}_t\colon H \, \to W$  by
$$ {R}_t \, = \, \frac{1}{2 \pi i }\, \int_{{\mathcal C}} \,
(A_{t} \,+ \, \lambda \, {\widetilde j}) ^{-1}\, d \lambda ,
$$we have  $P_t = {\widetilde j}\circ {R}_{t }$.  Moreover, if
$Q_{t } \, = \, R_{t } \, \circ \, {\widetilde j}, $  we have that
$Q^2_t= Q_t$ and hence each $Q_{t } $ is a projector belonging to
${\mathcal L }(W).$

By \cite[Chapter II, Section 6]{K} there exist two smooth paths $
U$ and $ V$ of unitary operators of $H$ and $W$ respectively,
defined in $ [t_j -\eta, t_j+\eta] $, such that $ U_{t_j} \, = \,
{\rm Id}_H$ ,  $V_{t_j} \, = \, {\rm Id}_W$ and such that
\begin{equation}\label{eq:conjugationjac1} P_{t } \, U_t \,
= \, U_t \, P_{t_j } ; \quad  Q_{t } \, V_t \, = \, V_t \, Q_{t_j
}. \end{equation}

Taking eventually a smaller $\eta,$ we can consider the smooth
operator valued function $ N_{z } \, = \, U_{t}^{-1}\, A_{z} \,
V_{t} $ defined on some  open neighborhood of the closed  domain
${\mathcal D}_ j = [t_j -\eta, t_j +\eta]\times[-1,1]$ together
with  the differential one-form
\begin{equation}\label{N}
          \theta \, = \, d N_{z }\, N_{z }^{-1}
\end{equation}

We claim that $ \theta $   take values in ${\mathcal T}(H),$ where
$\mathcal T(H)$  is the
  trace class,  and that
\begin{equation}\label{6.3}
{\rm Tr} \, d{\mathcal A}_z \, {\mathcal A}_z^{-1} \, ={\rm Tr}\,
d N_{z }\, N_{z }^{-1} .
\end{equation}

Indeed, denoting by dot the ordinary derivative with respect to
$t$ and using  ${\mathcal A}_z^{-1} = {\widetilde j} \circ
A_z^{-1}$ we obtain
\begin{eqnarray}\label{eq:formulonajac1}
d N_{z }\, N_{z }^{-1}  &=& -\, U_{t}^{-1} \, \dot{U}_{t} \, d t +
U_{t}^{-1} \, A_{z}\dot{V}_{t} \,V_{t}^{-1} \, A _{z }^{-1}\,U_{t}
d t\cr  &+& \, U_{t}^{-1} \, { \dot S}_t {\mathcal A}_{z}^{-1}\,
U_{t } dt + i U_{t}^{-1} {\mathcal A}_{z}^{-1}\, U_{t } ds .
\end{eqnarray}
The coefficients of all terms in the right hand side of
\equ{eq:formulonajac1} belong to the trace class. The last two
because ${\mathcal A}_{z}^{-1}\in {\mathcal T}(H).$ For the first
two let us recall that $U_{t}$ is a solution of the Cauchy problem
\begin{equation}\label{eq:commutatoreUjac1}
\left \{\begin{array}{ll}  \, \dot U_{t}  \, U_{t}^{-1}\, = \,
[\,  \dot P_{t} \, , P_{t}]\\
\noalign{\vskip 1 truemm} U_{t_j } \, = \, I
\end{array}\right.
\end{equation}
where $[\, \dot P_{t} \, , P_{t}]$ is the commutator. Since the
spectral subspace associated with the part of the spectrum of
${\mathcal A}_{t}$ lying in the interval $[- \mu , \mu ]$ is
finite dimensional it follows that  $[\, \dot P_{t } \, , P_{t }]$
is of trace class and indeed a finite rank operator. Thus by
\eqref{eq:commutatoreUjac1} we have that $\dot U_{t }=  [\, \dot
P_{t } \, , P_{t }]U_{t }$ is of  trace class and hence so is
$U_{t }^{-1}\dot U_{t }.$ The same argument shows that $
\dot{V}_{t} \, V_{t}^{-1}\in {\mathcal T}(H)$ and this completes
the proof of the first assertion.

In  order to show that \eqref{6.3} holds we first notice that by
definition of $R_t$
$$\dot{U}_{t} \,U_{t}^{-1}\, = \,[\dot{P}_{t} \, , \, P_{t }] \,
=\, {\tilde j} \, \dot{{R}}_{t }\, {\tilde j} \, {R}_{t }  \,-
  \, {\tilde j} \, {R}_{t} \, {\tilde j} \, \dot{{R}}_{t}\,= \,
  {\tilde j}D_t ,$$ where{ $D_{t} \, = \, \dot{R}_{t} \,
{\tilde j} \, R_{t}\, - \, R_{t} \, {\tilde j} \, \dot{R}_{t}. $
In the same way we get $\,\dot{V}_{t} V_{t}^{-1} \, =\, D_t{\tilde
j} .$

By commutativity of the trace
$${\rm Tr}\, U_{t }^{-1}\dot U_{t }= {\rm Tr} \,\dot{U}_{t} \,
U_{t}^{-1} \, = \,{\rm Tr} \,{\tilde j} \, D_{t} \, = \, {\rm Tr}
\,D_{t} \, {\tilde j}\,= \, {\rm Tr}\,\dot{V}_{t} \, V_{t}^{-1}.$$
 From this, taking in account that $ {\rm Tr}\, [ U_{t}^{-1} \,
A_{z} \, \dot{V}_{t} \, V_{t}^{-1} \, A_{z}^{-1}\,U_{t} ] \, = \,
{\rm Tr}\, \dot{V}_{t} \, V_{t}^{-1} $ it follows that
$${\rm Tr}[-\, U_{t}^{-1} \, \dot{U}_{t} \,  +  U_{t}^{-1} \, A_{z}
\dot{V}_{t} \,V_{t}^{-1} \, A _{z }^{-1}\,U_{t} ] \equiv 0.$$
Finally, substituting the above  in \equ{eq:formulonajac1} and
summing up, we obtain $$ {\rm Tr}\, d N_{z }\, N_{z }^{-1}  = {\rm
Tr} \, [U_{t}^{-1}( \, {\dot S}_t {\mathcal A}_{z}^{-1}\,  dt + i
{\mathcal A}_{z}^{-1}\, ds\, )U_{t }] = {\rm Tr} \, d {\mathcal
A}_{z} \, {\mathcal A}_{z}^{-1} ,$$ since the trace do not change
under conjugation. This proves \eqref{6.3}.

 From the very  definition of the one form $\theta \,=\, d N_{z }\,
N_{z }^{-1}$ it follows easily  that $d N_{z }\, N_{z }^{-1}
P_{t_j}\, = \, P_{t_j} d N_{z }\, N_{z }^{-1}$ for every $z\in
{\mathcal D}_j$ belonging to its domain and hence  $P_{t_j}$
reduces the coefficients of $\theta$.

Let $ H_{j} \, = \Imm\,P_{t_j} = \Ker \, {\mathcal A}_{t_j}.$
Under the splitting $ H\, =\, H_j \oplus H_j ^{\bot}$ we  can
write  $\theta$ as a  matrix of one forms
\begin{equation}\label{eq:6.7jac1}
\theta\, = \, \begin{pmatrix} \theta_0 & 0 \cr 0 & \theta_1
\end{pmatrix},
\end{equation}
where   $ \theta_0 = d N_{z } \, N_{z}^{-1}\, |_{H_j } $ and $
\theta_1 = d N_{z } \, N_{z}^{-1}\,|_{H_j ^{\bot}}.$

Taking traces,
\begin{equation}\label{eq:6.8jac1}
\frac{1}{2 \pi i } \,\int_{\partial {\mathcal D}_j} \,    {\rm Tr}
\, \theta \, = \, \frac{1}{2 \pi  i} \, \int_{\partial {\mathcal
D}_j} \,{\rm Tr}\, \theta_0 \,+\, \frac{1}{2 \pi i  } \,
\int_{\partial {\mathcal D}_j} \, {\rm Tr} \, \theta_1 .
\end{equation}

We claim that the last term in \eqref{eq:6.8jac1} vanishes. In
order to prove this  we  observe first that
${N_{z}^{-1}}|_{{H_{j}^{\bot}}}$ exists for all $z \, \in
{\mathcal D}_j$ and not only on the boundary.  On the other hand
since both $U_t\, {\tilde j}$ and ${\tilde j}\, V_t$ verify the
differential equation $\dot{W} = [\dot{P}_{t} \, , \, P_{t }] \,
W$ with the same initial condition $W(t_j)= {\tilde j}$ we have
that $U_t \,{\tilde j} ={\tilde j}\, V_t$ everywhere. It follows
then that $N_z=  U_{t}^{-1}A_t V_{t}\, +\, i\,s\, {\tilde j}.$
Substituting this into the expression $ dN_z
\,N_{z}^{-1}|_{H_{j}^{\bot}}$ and writing it in the form

$ E_z dt + F_z ds ,$ a  direct computation yields $
\partial_s E_z \, = \, \partial_t F_z.
$ Since ${\mathcal D}_j$ is simply connected  the one form $  {\rm
Tr}\, \theta_1$ is exact and  its integral over $\partial
{\mathcal D}_j$ vanishes.

Combining  with  \equ{A} and \eqref{6.3},  we obtain
\begin{equation} \label{integral}
\frac{1}{2 \pi  i} \, \int_{\partial {\mathcal D}_j} \,{\rm Tr}\,
d{{\mathcal A}}_{z} \, {\mathcal A}^{-1}_{z }\, = \,\frac{1}{2\pi
i } \int_{\partial {\mathcal D}_j} \,{\rm Tr}\, \left( dN_z
N_z^{-1} \right)\vert_{H_j},
\end{equation}
where the  right hand side is an integral of the trace form on a
finite dimensional space $H_j$.

Now let us turn our attention to the signature. At this point we
will identify the finite dimensional subspace $ H_j =
\Ker{\mathcal A}_{t_j}$ with $\Imm \, Q_{t_j}=\Ker\, N_{t_j}.$
With this identification   ${\tilde j}\colon \Imm \,Q_{t_j}\to
H_j$ becomes the identity. Let us consider the path of symmetric
endomorphisms $M  \colon [t_j -\eta, t_j + \eta ] \to {\mathcal L}
(H_j) $ given by $ M_t ={N_{t}}\vert_{H_j}$. Writing down the
definition of $N_t$ we find that $\forall u,v \in H_j$
\begin{eqnarray*}
\langle \dot{M}_t u,v\rangle_{L^2}&=& -\langle U_t^{-1} \dot{U}_t
U_t^{-1} A_t V_t u,v\rangle_{L^2}  +\\
&+& \langle U_t^{-1} \dot{A}_t  V_t u,v\rangle_{L^2} + \langle
U_t^{-1} A_t \dot{V}_tu,v\rangle_{L^2} .
\end{eqnarray*}
Putting $t=t_j$ in the above formula and using the fact that
${\mathcal A}_t$ is symmetric we notice that the first and the
last term in the right hand side of the above equation vanish.
Therefore
\begin{equation}\label{f2}
\sgn \,\Gamma(M_t, t_j) = \sgn \, \Gamma({\mathcal A}, t_j).
\end{equation}
In view of \eqref{f2} and \eqref{integral} Theorem \ref{TEOREMA}
follows from the following result

\begin{proposition}
If $M$  is the path defined above  and if $ M_z = M_t +is \Id,$
for $ z \in {\mathcal D}_j$, then
$$ \sgn \,\Gamma(M, t_j)   = \,\frac{1}{2\pi i} \int_{\partial
{\mathcal D}_j} \,{\rm Tr}\, dM_z M_z^{-1} .$$
\end{proposition}

\proof By   Kato's  Selection Theorem \cite[Chapter II, Theorem
6.8]{K} there exists smooth functions $\lambda_{1}(t) , \ldots ,
\lambda_{n_j}(t)$ representing for each $t$ the eigenvalues of the
symmetric matrix $M_t$. Equivalently, $M_t$ is similar to  a
smooth path of matrices $ \Delta_{t}$ having  the form $
\Delta_{t} \, = \,diag \, [ \lambda_{1}(t) \, , \ldots ,
\lambda_{n_j }(t) \,].$ By \cite[Chapter II, Theorem 5.4]{K}, $
{\rm Tr}\, \dot{M}_t = \sum_{i=1}^{n_j} \dot{\lambda_i}(t).$

  Putting $\Delta_{z}=
\Delta_{t} +is\,{\rm Id},$
\[ {\rm Tr}  \,dM_zM_z^{-1}=  \,{\rm Tr} \, d \Delta_{z} \,
\Delta _{z}^{-1} \,=  \sum _{l = 1}^{n_j} d\, ( \lambda_{l}(t) \,
+ i \,s)\, ( \lambda_{l}(t) \, + i \, s).^{-1}\] Now  by
elementary integration
$$\frac{1}{2 \pi  i}\int_{\partial {\mathcal D}_j}
d ( \lambda_{l}(t)  + i  s )\, ( \lambda_{l}(t)  + i  s)^{-1}=$$

\begin{equation}\label{fm}
= \left \{\begin{array}{lll}  \, \, \, \, \, \, \, 0 \, \, \, \, \, {\rm if}  & \mbox{$
\lambda_{l}(t_j \, - \, \eta )\,
\lambda_{l}(t_j \, + \, \eta )\, > \, 0$}\\
\noalign{\vskip 1 truemm}  \, \, \, \, \, \, \,1 \, \, \, \, \, {\rm if} & \mbox{$ \lambda_{l}(t_j \, -
\, \eta )\, \lambda_{l}(t_j \, + \, \eta )\,
  < \, 0 \,$ and $\, \lambda_{l}(t_j\, + \, \eta )\,> \,
0 $}\\ \noalign{\vskip 1 truemm}
  -1  \, \,  \,  \, {\rm if}  & \mbox{$
  \lambda_{l}(t_j \, - \, \eta )\, \lambda_{l}(t_j\, + \,
\eta )\, < \, 0 \, $ and $ \, \lambda_{l}(t_j\, + \, \eta )\,< \,
0 $}.
\end{array}\right.
\end{equation}
 Summing
over $l=1, \ldots , n_j$ in \equ{fm}  and using Proposition
\ref{prop5.1} we obtain
\begin{equation}
\frac{1}{2\pi i } \int_{\partial {\mathcal D}_j } {\rm Tr}\, dM_z
M_z^{-1} = \mu (M_{t_j -\eta} ) - \mu (M_{t_j +\eta} ) = \sgn
\,\Gamma(M, t_j).
\end{equation}
This completes the proof of the Theorem \ref{TEOREMA}.\qedhere

\section{Relation with the Maslov Index}\label{sec:maslovindex}

Let us consider again the path  $\tilde{{\mathcal A}}$ of real
unbounded  self-adjoint operators defined by \equ{AR}. According
to the results of  section \ref{sec:conjugateindex}, we have that
$t \in [0 , 1]$ is a conjugate instant along $\gamma$ if and only
if  $\ker \tilde{{\mathcal A}}_t\neq \{0\}.$ Let $\Psi_t (x)$ be
the solution to the initial value problem
\begin{equation} \label{flow}
\quad \ \left\{
\begin{array}{ll} \Psi_t'(x) =\sigma H_t (x) \Psi_t (x)  &\mbox{  }
x \in [0,1]\\ \Psi_t (0) = I,
\end{array} \right.
\end{equation}
where $\sigma $  and  $H_t (x)$ are matrices defined  as in
\equ{sigma} and \equ{acca} with $z=t, $ but  with real
coefficients this time.

The path $\Psi_t =\Psi_t (1) ,\, t\in[0 ,1]$ is a path of real
symplectic matrices i.e.\, $\Psi_t \sigma \Psi_t^*= \sigma$. Since
the symplectic  group acts on the manifold $\Lambda_n$ of all
Lagrangian subspaces of $\R^{2n},$ the action of the path $\Psi $
on fixed Lagrangian subspace $l$ produces a path on $\Lambda_n$.
This path can be used in order to count conjugate points of
$\gamma$.

Let us recall that in terms of the complex structure $\sigma$ the
manifold $\Lambda_n $ can be defined as the submanifold  of the
Grassmanian $G_n (2n) $ whose elements are $n$-dimensional
subspaces $l$ of $\R^{2n}$ such that $\sigma( l) $ is orthogonal
to $l.$ Given a path of Lagrangian subspaces $ \lambda \colon
[a,b] \to \Lambda_n $ and a fixed Lagrangian subspace $l$  such
that $\lambda _a$ and $\lambda _b$ are transverse to $l$,  there
is a well defined integral-valued homotopy invariant
$\mu_l\,(\lambda, [a,b]), $ called {\em Maslov index\/}
\cite{RS2},  which counts algebraically the number of crossing
points of $\lambda$ with $l$ i.e.\, points $t\in [a,b]$ at which
$\lambda_t  $ fails to be transverse  to $l$. Here we take $ l =
\{ 0\} \times \R^n$ and $\lambda_t = \Psi_t ( l)$.  With this
choice, we have that  $t $ is a conjugate instant along $\gamma$
if and only if $\lambda_t \cap l \neq \{0\}$. Since conjugate
instants cannot accumulate at $0$ we can find an $\ve > 0$ such
that there are no conjugate instants in $[0,\ve]$. The Maslov
index  of $\lambda \colon  [\ve ,  1] \to \Lambda_n $ is well
defined and independent of the choice of $\ve.$

The {\em Maslov index\/} of a $p$-geodesic is defined by
$\iMaslov(\gamma) = \mu_l\,(\lambda, [\ve ,1]).$
\begin{proposition}
\begin{equation}\label{mas=spec}
\isf(\gamma) = \icon(\gamma)=\iMaslov(\gamma).
\end{equation}
\end{proposition}

\proof   We will be rather sketchy here since the arguments are
very close to those used  in the first part of the previous
theorem.  A different approach in the geodesic case  can  be  found
in \cite{PPT}.

There is a  natural identification of  the  tangent space  to the
manifold $\Lambda (n)$ at a given point $l$ in $\Lambda (n)$  with
the set of quadratic  forms on $l$ which allows to express the
Maslov index of a generic path as a sum  of  the signatures of the
crossing forms at points of regular crossing  analogous to that in
Proposition \ref{crossform} for the spectral flow.  In terms of
this identification the crossing  form $\Gamma(\lambda, t)$  of a
path $\lambda$ at a given point $t$ can be constructed  as
follows: let $ M \colon U \to \mathcal{L}(\RR^n;\RR^{2n})$  be a
smooth path of monomorphisms defined on  a neighborhood  $U$ of
$t$  such that ${\rm Im} \, M_s = \lambda_s.$  By  \cite{RS2}, the
crossing form at $t$ is the quadratic form defined by

\begin{equation}\label{mascr}
\Gamma(\lambda, t) v= \langle \sigma\dot{M}_t w, M_tw \rangle \
\hbox { for } \,v\in \lambda_t \cap l;\, M_t w =v .
\end{equation}
The resulting  form is independent of the choice of  the frame
$M$.

As before, a  crossing point  $t$ is regular if the crossing form
$\Gamma(\lambda, t)$ is non-degenerate.  For paths with only
regular crossing  points the Maslov index can be computed by
$\mu_{l}(\lambda, [a,b]) = \sum_t \Gamma(\lambda, t).$

From the above discussion  the paths $\hat{h}_t, \,
\tilde{{\mathcal A}}_t$  and the path $\lambda_t = \Psi_t ( l)$
defined in $[\ve ,1]$ have the same crossing points. The crossing
form $\Gamma(\lambda, t)$ can be easily computed in this case
taking as $M$  the second column of the block decomposition of
$$\Psi_t \,= \,  \begin{pmatrix} a_t &
b_t \cr c_t & d_t
\end{pmatrix}.  $$
Then $ M_t  w \in l$ if and only if  $b_t w =0.$

It turns out that $\lambda_t \cap l =  \{ (0,v )|  v = d_t w \}\,
\hbox{\rm and}$
$$\Gamma(\lambda, t)(0, v)  = - \langle\dot{b}_t w,
d_t w \rangle  =  -  \langle Jd_t w,d_t w\rangle  = -  \langle J
v,v\rangle $$ by \eqref{flow}. Identifying $\lambda_t \cap l$ with
$V_t= \Imm \, d_t,$ the crossing form is $-\langle J v,v\rangle$.

On the other hand, by \eqref{eq:crossjac1}, for   $ u \in \ker
\mathcal A_t$
\begin{equation}
\Gamma( {\mathcal A }, t) u= - \int_0^1 \langle{\dot S}_{t}(x)
u(x) , u(x)\rangle dx.
\end{equation}
where as before $\cdot{ } $ denotes the derivative with respect to
$t$.

For any $s\in(0,1]$, the function $u_s (x)=u(\frac st \,x)$ solves
the Cauchy problem
\begin{equation}\label{cp}
\left\{ \begin{array}{ll} Ju_s ''(x) + S_s (x)u_s (x)
= 0&\mbox{ for all } x \in [0,1] \\
u_s (0)=0 , \quad u_s' (0) = \frac st u' (0)
\end{array}\right.
\end{equation}

If we differentiate the equation in \eqref{cp} with respect to $s$
and we evaluate at $s=t$, we get

\begin{equation*}
\left\{\begin{array}{ll} J{\dot u}_{t} ''(x) + {\dot S}_{t}
(x)u_{t} (x) + S_{t} (x){\dot
u}_{t} (x)= 0 \quad \mbox{ for all } s \in [0,1]\\
{\dot u}_{t} (0)=0 , \quad {\dot u}_{t}' (0) = 0.
\end{array}\right.
\end{equation*}
Taking into account that $ u \in \ker \widetilde{\mathcal A}_t$
and the previous computation, integrating by parts we  get
\begin{eqnarray*}
\Gamma(\widetilde{\mathcal A} , t) (u)&=& \int_0^1 \langle J{\dot
u}_{t}'' + S_{t} {\dot u}_t(x) , u_t (x)\rangle dx=- \langle J
u_t'(1) , {\dot u_t} (1) \rangle
\end{eqnarray*}

Since  $u_t(x)=u(x), \, \dot{u}_t(x) = \frac xt u'(x) $ we get
$$ \Gamma( \tilde{\mathcal A }, t) (u) = -
\frac1t  \langle J u' (1) , u' (1) \rangle.$$

The above calculation shows  that the isomorphism sending  $ u \in
\ker \tilde{{\mathcal A}}_t $ into  $t u'(1) \in V_t $ transforms
$\Gamma( \tilde{{\mathcal A}}, t)$ into $\Gamma(\lambda, t).$ We
conclude that the regular crossings of $\tilde{{\mathcal A }}$
correspond to the regular crossings of $\lambda$ and moreover the
crossing forms have the same signature. Now the assertion
\equ{mas=spec} follows from \equ{mascr} and the perturbation
argument used in the proof of the main theorem. \qedhere
\smallskip

We close this section by showing, as promised, that in the case of
a regular geodesic $\icon(\gamma) = \ispec(\gamma)$ coincides with
the expression given by formula \equ{icon}.

Let us observe  that  if $\xi\in {\mathcal I}$ then for any Jacobi
field  $\eta$ along $\gamma$ the function given by
$g(\frac{D}{dx}\xi(x),\eta(x))$ is constant (to prove this it is
enough to use \equ{2.5} in the expression for the derivative of
this function). In particular, if also $\eta\in {\mathcal I}$ then
$g(\frac{D}{dx}\xi(x),\eta(x)) \equiv 0$. This shows that for any
Jacobi field $\xi \in \ker \tilde{{\mathcal A}}_t $ the covariant derivative $\Ddt\xi(t)$ belongs
to  ${\mathcal I}[t]^{\bot}.$

Let $M\colon \ker\tilde{{\mathcal A}}_t \to {\mathcal
I}[t]^{\bot}$ be the monomorphism defined by
$$M u = t  \Ddt \big(\sum_{i=1}^{n} u_{i}(\frac xt)  e^
{i}(x)\big).$$ By  dimension counting $M$ is an isomorphism. On
the other hand, using our previous calculation  we get
$$ g(Mu, Mu) = -\langle Ju'(1),u'(1)\rangle =  t
\Gamma(\tilde{{\mathcal A}}, t).$$

Taking signatures and summing over all conjugate points we obtain
the desired conclusion.

\end{document}